
\input amstex
\documentstyle{amsppt}
\pageheight{56pc}
\pagewidth{38pc}
    \def\nz{\hbox{\text{0}{ \raise 3pt \hbox{\kern -10pt
    \vrule width 8pt height 0.4 pt }\kern 2pt}}}

\topmatter
\centerline{{\bf JACOBIAN OF MEROMORPHIC CURVES}
\footnote"*"{
Abhyankar's work was partly supported by NSF Grant DMS 91-01424 and
NSA grant MDA 904-97-1-0010.}}
\centerline{By}
\centerline{Shreeram S. Abhyankar and Abdallah Assi}
\centerline{Mathematics Department, Purdue University, 
West Lafayette, IN 47907, USA;}
\centerline{e-mail: ram\@cs.purdue.edu}
\centerline{Universit\'e d'Angers, Math\'ematiques, 49045 Angers cedex 01,
France;}
\centerline{e-mail: assi\@univ-angers.fr}
\abstract{The contact structure of two meromorphic curves gives a 
factorization of their jacobian.}
\endabstract
\endtopmatter

\document


\centerline{{\bf Section 1: Introduction}}

Let $J(F,G)=J_{(X,Y)}(F,G)$ be the jacobian of $F=F(X,Y)$ and
$G=G(X,Y)$ with respect to $X$ and $Y$, i.e., let
$J(F,G)=F_XG_Y-F_YG_X$ where subscripts denote partial derivatives.
Here, to begin with, $F$ and $G$ are plane curves, i.e., polynomials
in $X$ and $Y$ over an algebraically closed ground field $k$ of
characteristic zero. More generally, we let $F$ and $G$ be meromorphic
curves, i.e., polynomials in $Y$ over the (formal) meromorphic series
field $k((X))$.

In terms of the the contact structure of $F$ and $G$, we shall produce
a factorization of $J(F,G)$. 
Note that if $G=-X$ then $J(F,G)=F_Y$; in this special case, our results 
generalize some results of Merle \cite{Me}, Delgado \cite{De}, and Kuo-Lu 
\cite{KL} who studied the situation when $F$ has one (Merle) or two (Delgado) 
or more (Kuo-Lu) branches. These authors restricted their attention to
the analytic case, i.e., when $F$ is a polynomial in $Y$ over the (formal) 
power series ring $k[[X]]$. With an eye on the Jacobian Conjecture, we
are particularly interested in the meromorphic case.

The main technique we use is the method of {\bf Newton Polygon}, i.e.,
the method of deformations, characteristic sequences, truncations, and 
contact sets given in Abhyankar's 1977 Kyoto paper \cite{Ab}. 
In Sections 2 to 5 we shall review the relevant material
from \cite{Ab}. In Section 6 we shall introduce the tree of contacts and 
in Sections 7 to 9 we shall show how this gives rise to the factorizations.

The said Jacobian Conjecture predicts that if the Jacobian of two 
bivariate polynomials $F(X,Y)$ and $G(X,Y)$ is a nonzero constant then
the variables $X$ and $Y$ can be expressed as polynomials in $F$ and $G$,
i.e., if $0\ne J(F,G)\in k$ for $F$ and $G$ in $k[X,Y]$ then
$k[F,G]=k[X,Y]$. We hope that the results of this paper may contribute
towards a better understanding of this Bivariate Conjecture, and hence
also of its obvious Multivariate Incarnation.

\centerline{}

\centerline{{\bf Section 2: Deformations}}

We are interested in studying polynomials in indeterminates $X$ and $Y$ over
an algebraically closed ground {\bf field } $k$ of characteristic zero. 
To have more elbow room to maneuver, we consider the larger {\bf ring} 
$R=k((X))[Y]$ of polynomials in $Y$ over $k((X))$, i.e., with coefficients in 
$k((X))$, where $k((X))$ is the meromorphic series field in $X$ over $k$.

Given any
$$
g=g(X,Y)=\sum_{i\in\Bbb Z}g^{[i]}X^i=\sum_{j\in\Bbb Z}g^{((j))}Y^j\in R
$$
with
$$
g^{[i]}=g^{[i]}(Y)\in k[Y]
\quad\text{ and }\quad
g^{((j))}=g^{((j))}(X)\in k((X))
$$
we put
$$
\text{Supp}_Xg=\{i\in\Bbb Z:g^{[i]}\ne 0\}
\quad\text{ and }\quad
\text{Supp}_Yg=\{j\in\Bbb Z:g^{((j))}\ne 0\}
$$
and we call this the $X${\bf -support} and the $Y${\bf -support} of $g$ 
respectively. Note that these supports are bounded from below and above 
respectively, and upon letting
$$
\gamma^{\sharp}=\text{ord}_Xg=\text{the $X${\bf -order} of }g
\quad\text{ and }\quad
\gamma=\text{deg}_Yg=\text{the $Y${\bf -degree} of }g
$$
we have
$$
\gamma^{\sharp}=\cases
\text{min}(\text{Supp}_Xg)&\text{if }g\ne 0\\
\infty&\text{if }g=0
\endcases
\quad\text{ and }\quad
\gamma=\cases
\text{max}(\text{Supp}_Yg)&\text{if }g\ne 0\\
-\infty&\text{if }g=0.
\endcases
$$
Now
$$
g^{[i]}(Y)=\sum_{j\in\Bbb Z}g^{(i,j)}Y^j
\quad\text{ and }\quad
g^{((j))}(X)=\sum_{i\in\Bbb Z}g^{(i,j)}X^i
\quad\text{ with }\quad
g^{(i,j)}\in k
$$
and we put
$$
\text{Supp}(g)=\text{Supp}_{(X,Y)}g=\{(i,j)\in\Bbb Z\times\Bbb Z:
g^{(i,j)}\ne 0\}
$$
and we call this the {\bf support}, or the $(X,Y)${\bf -support}, of $g$.
We put
$$
\text{inco}_Xg=\cases
g^{[\gamma^{\sharp}]}&\text{if }g\ne 0\\
0&\text{if }g=0
\endcases
\quad\text{ and }\quad
\text{deco}_Yg=\cases
g^{((\gamma))}&\text{if }g\ne 0\\
0&\text{if }g=0
\endcases
$$
and we call this the {\bf $X$-initial-coefficient} and the
{\bf $Y$-degree-coefficient} of $g$ respectively.
Upon letting
$$
\align
\widehat\gamma^{\sharp}&=\text{ord}(g)=\text{the (total) {\bf order} of }g\\
&=\text{ord}_{(X,Y)}g=\text{the $(X,Y)${\bf -order} of }g
\endalign
$$
we have
$$
\widehat\gamma^{\sharp}=\cases
\text{min}\{i+j:(i,j)\in\text{Supp}(g)\}&\text{if }g\ne 0\\
\infty&\text{if }g=0
\endcases
$$
and we put
$$
\text{info}(g)=\text{info}_{(X,Y)}g=\cases
\sum_{i+j=\widehat\gamma^{\sharp}}g^{(i,j)}X^iY^j&\text{if }g\ne 0\\
0&\text{if }g=0
\endcases
$$
and we call this the {\bf initial-form}, or the {\bf $(X,Y)$-initial-form}, 
of $g$. 
\footnote{In an obvious manner, the definitions of 
Supp$_Xg$, ord$_Xg$, inco$_Xg$, Supp$_Yg$, deg$_Yg$, deco$_Yg$, 
Supp$_{(X,Y)}g$, ord$_{(X,Y)}g$, and info$_{(X,Y)}g$,
can be extended to any $g$ in $k((X))[Y,Y^{-1}]$, and for any such $g$ we can
also define ord$_Yg$ and inco$_Yg$, and then we have:
$g=0\Leftrightarrow\text{Supp}_Xg=\emptyset\Leftrightarrow\text{ord}_Xg=\infty
\Leftrightarrow\text{inco}_Xg=0\Leftrightarrow\text{Supp}_Yg=\emptyset
\Leftrightarrow\text{ord}_Yg=\infty\Leftrightarrow\text{deg}_Yg=-\infty
\Leftrightarrow\text{inco}_Yg=0\Leftrightarrow\text{deco}_Yg=0
\Leftrightarrow\text{Supp}_{(X,Y)}g=\emptyset
\Leftrightarrow\text{ord}_{(X,Y)}g=\infty
\Leftrightarrow\text{info}_{(X,Y)}g=0$.}
If $g\in k[X,Y]$ then upon letting
$$
\align
\widehat\gamma&=\text{deg}(g)=\text{the (total) {\bf degree} of }g\\
&=\text{deg}_{(X,Y)}g=\text{the $(X,Y)${\bf -degree} of }g
\endalign
$$
we have
$$
\widehat\gamma=\cases
\text{max}\{i+j:(i,j)\in\text{Supp}(g)\}&\text{if }g\ne 0\\
-\infty&\text{if }g=0
\endcases
$$
and we put
$$
\text{defo}(g)=\text{defo}_{(X,Y)}g=\cases
\sum_{i+j=\widehat\gamma}g^{(i,j)}X^iY^j&\text{if }g\ne 0\\
0&\text{if }g=0
\endcases
$$
and we call this the {\bf degree-form}, or the {\bf $(X,Y)$-degree-form}, 
of $g$. 
\footnote{Again, in an obvious manner, the definitions of 
deg$_{(X,Y)}g$ and defo$_{(X,Y)}g$ can be extended to any $g$ in 
$k[X,X^{-1},Y,Y^{-1}]$, and for any such $g$ we can
also define deg$_Xg$ and deco$_Xg$, and then we have:
$g=0\Leftrightarrow\text{deg}_Xg=-\infty\Leftrightarrow\text{deco}_Xg=0
\Leftrightarrow\text{deg}_{(X,Y)}g=-\infty
\Leftrightarrow\text{defo}_{(X,Y)}g=0$.}

Given any $z=z(X)\in k((X))$, we write
$$
z=z(X)=\sum_{i\in\Bbb Z}z[i]X^i
\quad\text{ with }\quad z[i]\in k
$$
and we put
$$
z[i]=0\text{ for all }i\in\Bbb Q\setminus\Bbb Z
$$
and we let
$$
\epsilon(z)=\cases
\text{the set of all $(U,V,W)\in\Bbb Z^3$ such that $U>0<V$}\\
\text{ and $iV/U\in\Bbb Z$ for all $i\in\text{Supp}_Xz$ with $i<WU/V$}
\endcases
$$
and we call this the {\bf edge} of $z$, 
and for any $(U,V,W)\in\epsilon(z)$ we let
$$
z^{\dagger}(X,U,V,W)=\sum_{i<WU/V}z[i]X^{iV/U}\in k((X))
$$
and
$$
z^{\ddagger}(X,U,V,W,Y)=z^{\dagger}(X,U,V,W)+X^W Y\in R
$$
and
$$
z^{\dagger *}(X,U,V,W)=\sum_{i\le WU/V}z[i]X^{iV/U}\in k((X))
$$
and we call these the {\bf $(U,V,W)$-truncation},
the {\bf $(U,V,W)$-deformation}, and the
{\bf strict $(U,V,W)$-truncation} of $z(X)$ respectively.
Given any $H=H(X,Y)\in R$, we are interested in calculating
ord$_XH(X^V,z^{\ddagger}(X,U,V,W,Y))$ and 
inco$_XH(X^V,z^{\ddagger}(X,U,V,W,Y))$. 
\footnote{To motivate the definitions of $\epsilon(z)$ and $z^{\ddagger}$, 
given any $H=H(X,Y)=\sum_{i,j}H^{(i,j)}X^iY^j\in R$ with $H^{(i,j)}\in k$, 
let $\Gamma^{\sharp}$ and $\Theta^{\sharp}(Y)$ be the values of 
ord$_XH(X^V,z^{\ddagger}(X,U,V,W,Y))$ and 
inco$_XH(X^V,z^{\ddagger}(X,U,V,W,Y))$ when $z=0$, i.e., let
$\Gamma^{\sharp}=\text{ord}_XH(X^V,X^WY)$ and 
$\Theta^{\sharp}(Y)=\text{inco}_XH(X^V,X^WY)$.
Also let $\Gamma$ and $\Theta(X,Y)$ be the
weighted order and the weighted initial form of $H(X,Y)$, 
when we give weights $(V,W)$ to $(X,Y)$, i.e., let
$\Gamma=\text{min}\{iV+jW:(i,j)\in\text{Supp}_{(X,Y)}H(X,Y)\}$ and 
$\Theta(X,Y)=\sum_{iV+jW=\Gamma}H^{(i,j)}X^iY^j$. Then 
$\Gamma^{\sharp}=\Gamma=\text{ord}_{(X,Y)}H(X^V,Y^W)$ and 
$\Theta(X,Y)=\text{info}_{(X,Y)}H(X^V,Y^W)$. 
Moreover, $\Theta^{\sharp}(Y)$ and 
$\Theta(X,Y)$ determine each other by the formulas 
$\Theta^{\sharp}(Y)=\Theta(1,Y)$ and 
$\Theta(X,Y)=X^{\Gamma/V}\Theta^{\sharp}(X^{-W/V}Y)$.
The parameter $U$ is a normalizing parameter which essentially says that
we want to intersect the ``meromorphic curve'' $H(X,Y)=0$ with a deformation 
of the ``irreducible meromorphic curve'' $f(X,Y)=0$ where $f(X,Y)$ is a monic
irreducible polynomial of degree $U=n$ in $Y$ over $k((X))$; to do this we 
take a ``fractional meromorphic'' root $y(X)$ of $f(X,Y)=0$ with
$y(X^n)=z(X)\in k((X))$, and then after ``deforming'' $y(X)$ at 
$X^{W/U}$ we substitute the deformation in $f(X,Y)$ for $Y$. 
For further motivation see the definitions
of $\epsilon(f,\lambda)$ and $t(f,\lambda)$ displayed in the
middle of the next Section.} 
For this purpose we proceed to give a review on characteristic sequences.

\centerline{}

\centerline{{\bf Section 3: Characteristic Sequences}}

Let $\widehat R^{\natural}$ be the set of all {\bf monic polynomials} in $Y$ 
over $k((X))$, i.e., those nonzero members of $R$ in whom the coefficient of 
the highest $Y$-degree term is $1$. Let $R^{\natural}$ be 
the set of all {\bf irreducible monic polynomials} in $Y$ over $k((X))$, i.e., 
those members of $\widehat R^{\natural}$ which generate prime ideals in $R$; 
note that their $Y$-degrees are positive integers.

Given any $f=f(X,Y)\in R^{\natural}$ of $Y$-degree $n$, 
by Newton's Theorem
$$
f(X^n,Y)=\prod_{1\le j\le n}[Y-z_j(X)]
\quad\text{ with }\quad z_j(X)\in k((X))
$$
where we note that Supp$_Xz_j$ is independent of $j$. Let 
$m(f)=m_i(f)_{0\le i\le h(m(f))+1}$ be 
the {\bf newtonian sequence of characteristic exponents} of $f$ 
relative to $n$ as defined on page 300 of \cite{Ab}, let 
$d(m(f))=d_i(m(f))_{0\le i\le h(d(m(f)))+2}$
be the {\bf GCD-sequence} of $m(f)$ as defined on page 297 of \cite{Ab}, 
let $q(m(f))=q_i(m(f))_{0\le i\le h(q(m(f)))+1}$ be the 
{\bf difference sequence} of $m(f)$ as defined on page 301 of \cite{Ab}, let 
$s(q(m(f)))=s_i(q(m(f)))_{0\le i\le h(s(q(m(f))))+1}$ be the {\bf inner 
product sequence} of $q(m(f))$ as defined on page 302 of \cite{Ab}, and let 
$r(q(m(f)))=r_i(q(m(f)))_{0\le i\le h(r(q(m(f))))+1}$ be the {\bf normalized 
inner product sequence} of $q(m(f))$ as defined on page 302 of \cite{Ab}.
\footnote{It is really not necessary to look up \cite{Ab} for the
definitions of the sequences $m,d,q,s,r$, since they are completely
redefined in the next three sentences ending with the displayed
item $(\bullet)$.}
Note that then 
$$
\align
h(d(m(f)))&=h(m(f))=h(q(m(f)))\\
&=h(s(q(m(f))))=h(r(q(m(f))))=\text{a nonnegative integer}
\endalign
$$
and 
$$
d_0(m(f))=0
\quad\text{ and }\quad
d_{h(d(m(f)))+1}(m(f))=1
$$
and
$$
d_1(m(f))=m_0(f)=q_0(m(f))=s_0(q(m(f)))=r_0(q(m(f)))=n
$$
and
$$
\align
d_{h(m(f))+2}(m(f))&=m_{h(m(f))+1}(f)=q_{h(m(f))+1}(m(f))\\
&=s_{h(m(f))+1}(q(m(f)))=r_{h(m(f))+1}(q(m(f)))=\infty
\endalign
$$
and
$$
m_1(f)=q_1(m(f))=s_1(q(m(f)))/n=r_1(q(m(f)))=\text{min}(\text{Supp}_Xz_1)
$$ 
with the understanding that the min of the empty set is $\infty$. 
Also note that: $h(m(f))=0\Leftrightarrow f(X,Y)=Y$.
Finally note that if $f(X,Y)\ne Y$ then for $2\le i\le h(m(f))$ we have that
$d_i(m(f))$, $m_i(f)$, $q_i(m(f))$, $s_i(q(m(f)))$, $r_i(q(m(f)))$ are 
integers with $d_i(m(f))>0$ such that: 
$$
\cases
d_i(m(f))=\text{GCD}(m_0(f),m_1(f),\dots,m_{i-1}(f)),\\
m_i(f)=\text{min}(\text{Supp}_Xz_1\setminus d_i(m(f))\Bbb Z),\\
q_i(m(f))=m_i(f)-m_{i-1}(f),\\ 
s_i(q(m(f)))=q_1(m(f))d_1(m(f))+\dots+q_i(m(f))d_i(m(f)),\\
\text{and }r_i(q(m(f)))=s_i(q(m(f)))/d_i(m(f)).
\endcases
\tag$\bullet$
$$
In the rest of this Section we shall use the abbreviations
$$
d_i=d_i(m(f))\quad\text{ and }\quad s_i=s_i(q(m(f)))
$$
for all relevant values of $i$.
Let the sequence
$c(f)=c_i(f)_{1\le i\le h(c(f))}$ be defined by putting
$$
h(c(f))=h(m(f))\quad\text{ and }\quad
c_i(f)=m_i(f)/n\text{ for }1\le i\le h(c(f))
$$ 
and let us call this the {\bf normalized characteristic sequence} of 
$f$. Note that then 
$c_1(f)<c_2(f)<\dots<c_{h(c(f))}(f)$ are rational 
numbers, out of which only $c_1(f)$ could be an integer. To obtain an
alternative characterization of the noninteger members of this sequence, 
{\bf for any rational number} $\lambda$, we let
$$
p(f,\lambda)=\cases
\text{the unique nonnegative integer $\le h(c(f))$ such that}\\
c_i(f)<\lambda\le c_{j}(f)\text{ for }
1\le i\le p(f,\lambda)<j\le h(c(f))
\endcases
$$
and
$$
p^*(f,\lambda)=\cases
\text{the unique nonnegative integer $\le h(c(f))$ such that}\\
c_i(f)\le \lambda<c_{j}(f)\text{ for }
1\le i\le p^*(f,\lambda)<j\le h(c(f))
\endcases
$$
and
$$
D(f,\lambda)=n/d_{p+1}\quad\text{ with }\quad p=p(f,\lambda)
$$
and
$$
D^*(f,\lambda)=n/d_{p^*+1}\quad\text{ with }\quad p^*=p^*(f,\lambda)
$$
and
$$
S(f,\lambda)=\cases
(s_p+(n\lambda-m_p(f))d_{p+1})/n^2
&\text{if }p=p(f,\lambda)\ne 0\\
\lambda&\text{if }p=p(f,\lambda)=0
\endcases
$$
and, for any $z\in k((X))$, we let
$$
A(f,\lambda,z)=\prod_{i=1}^{p(f,\lambda)}\left((d_i/d_{i+1})
z[m_i(f)]^{(d_i/d_{i+1})-1}\right)^{d_{i+1}D(f,\lambda)/n}
$$
and
$$
\widehat A(f,\lambda,z)=A(f,\lambda,z)^{n/D(f,\lambda)}
$$
and
$$
E(f,\lambda,z,Y)=Y^{D^*(f,\lambda)/D(f,\lambda)}
-z[n\lambda]^{D^*(f,\lambda)/D(f,\lambda)}
$$
and
$$
\widehat E(f,\lambda,z,Y)=E(f,\lambda,z,Y)^{n/D^*(f,\lambda)}
$$
and we call these the {\bf $\lambda$-position}, the 
{\bf strict $\lambda$-position}, the {\bf $\lambda$-degree}, the 
{\bf strict $\lambda$-degree}, the {\bf $\lambda$-strength}, 
the {\bf $(\lambda,z)$-reduced-constant}, the 
{\bf $(\lambda,z)$-constant}, the {\bf $(\lambda,z)$-reduced-polynomial}, 
and the {\bf $(\lambda,z)$-polynomial} of $f$ respectively; note that the 
above objects $p$, $p^*$, $D$, $D^*$, $S$, $A$, 
$\widehat A$, $E$, and $\widehat E$ respectively correspond to the objects
$p(<)$, $p(\le)$, $D$, $E$, $s$, $B$, $\widehat B$, $P$, and 
$\widehat P$ introduced on pages 326-328 of \cite{Ab}. 
We also define the sequence
$m(f,\lambda)=m_i(f,\lambda)_{0\le i\le h(m(f,\lambda))+1}$ by putting
$$
h(m(f,\lambda))=p(f,\lambda)\quad\text{ and }\quad
m_i(f,\lambda)=m_i(f)D(f,\lambda)/n
\text{ for }0\le i\le p(f,\lambda)+1
$$
with the understanding that $m_i(f,\lambda)=\infty$ for 
$i=p(f,\lambda)+1$, and we define the sequence
$m^*(f,\lambda)=m^*_i(f,\lambda)_{0\le i\le h(m^*(f,\lambda))+1}$ by putting
$$
h(m^*(f,\lambda))=p^*(f,\lambda)\quad\text{ and }\quad
m^*_i(f,\lambda)=m_i(f)D^*(f,\lambda)/n
\text{ for }0\le i\le p^*(f,\lambda)+1
$$
with the understanding that $m^*_i(f,\lambda)=\infty$ for 
$i=p^*(f,\lambda)+1$, 
and we define the sequence 
$c(f,\lambda)=c_i(f,\lambda)_{1\le i\le h(c(f,\lambda))}$ by putting
$$
h(c(f,\lambda))=p(f,\lambda)\quad\text{ and }\quad
c_i(f,\lambda)=c_i(f)\text{ for }1\le i\le p(f,\lambda)
$$ 
and we define the sequence 
$c^*(f,\lambda)=c^*_i(f,\lambda)_{1\le i\le h(c^*(f,\lambda))}$ by putting
$$
h(c^*(f,\lambda))=p^*(f,\lambda)\quad\text{ and }\quad
c^*_i(f,\lambda)=c_i(f)\text{ for }1\le i\le p^*(f,\lambda)
$$ 
and we call these sequences
the {\bf $\lambda$-characteristic-sequence},
the {\bf strict $\lambda$-characteristic-sequence},
the {\bf $\lambda$-normalized-characteristic-sequence},
and the {\bf strict $\lambda$-normalized-characteristic-sequence}
of $f$ respectively.
We also let
$$
\epsilon(f,\lambda)=\cases
\text{the set of all }(z,U,V,W)\in k((X))\times\Bbb Z^3
\text{ such that $U=n$, $W/V=\lambda$,}\\
\text{and $(U,V,W)\in\epsilon(z)$ where $z=z(X)\in k((X))$ with }
f(X^n,z(X))=0
\endcases
$$
and we call this the {\bf $\lambda$-edge} of $f$.
Finally we define 
$t(f,\lambda)=t(f,\lambda)(X,Y)$ to be the unique
member of $R^{\natural}$ such that
$$
t(f,\lambda)(X^V,z^{\dagger}(X,U,V,W))=0
\text{ for some (and hence for all) }(z,U,V,W)\in\epsilon(f,\lambda)
$$
and we call this the {\bf $\lambda$-normalized-truncation} of $f$, and
we define 
$t^*(f,\lambda)=t^*(f,\lambda)(X,Y)$ to be the unique
member of $R^{\natural}$ such that
$$
t^*(f,\lambda)(X^V,z^{\dagger *}(X,U,V,W))=0
\text{ for some (and hence for all) }(z,U,V,W)\in\epsilon(f,\lambda)
$$
and we call this the {\bf strict $\lambda$-normalized-truncation} of $f$;
note that on page 294 of \cite{Ab} we have called these the open and closed
$(n\lambda)$-truncations of $f$ respectively. 

>From the above definitions of the various objects, we see that
$$
\cases
\text{$p(f,\lambda)$ and $p^*(f,\lambda)$ are integers with}\\
\text{$0\le p(f,\lambda)\le p^*(f,\lambda)\le h(m(f))$, and}\\
\text{$D(f,\lambda)$ and $D^*(f,\lambda)$ are positive integers with}\\
\text{$n/D^*(f,\lambda)\in\Bbb Z$ and 
$D^*(f,\lambda)/D(f,\lambda)\in\Bbb Z$}\\
\endcases
\tag{NP1}
$$
and
$$
\cases
\text{$t(f,\lambda)$ and $t^*(f,\lambda)$ are elements of 
$R^{\natural}$ such that:}\\
\text{$m(t(f,\lambda))=m(f,\lambda)$ and
$m(t^*(f,\lambda))=m^*(f,\lambda)$,}\\
\text{$c(t(f,\lambda))=c(f,\lambda)$ and
$c(t^*(f,\lambda))=c^*(f,\lambda)$,}\\
\text{deg$_Yt(f,\lambda)=D(f,\lambda)$ and
deg$_Yt^*(f,\lambda)=D^*(f,\lambda)$, and}\\
\text{$h(m(t(f,\lambda)))=h(c(t(f,\lambda)))=p(f,\lambda)$ and
$h(m(t^*(f,\lambda)))=h(c(t^*(f,\lambda)))=p^*(f,\lambda)$}
\endcases
\tag{NP2}
$$
and
$$
\cases
A(f,\lambda,z)\in k\text{ and }
\widehat A(f,\lambda,z)=A(f,\lambda,z)^{n/D(f,\lambda)}\in k\text{ are
such that:}\\
\text{if }f(X^n,z(X))=0\text{ then } 
A(f,\lambda,z)\ne 0\ne\widehat A(f,\lambda,z)\\
\endcases
\tag{NP3}
$$
and
$$
\cases
E(f,\lambda,z,Y)=Y^{D^*(f,\lambda)/D(f,\lambda)}
-z[n\lambda]^{D^*(f,\lambda)/D(f,\lambda)}\in k[Y]\\
\text{and }\widehat E(f,\lambda,z,Y)
=E(f,\lambda,z,Y)^{n/D^*(f,\lambda)}\in k[Y]\\
\text{are monic polynomials of degrees
$D^*(f,\lambda)/D(f,\lambda)$ and
$n/D(f,\lambda)$ respectively,}\\
\text{where }z[n\lambda]\in k\text{ is such that: }
z[n\lambda]\ne 0\Leftrightarrow n\lambda\in\text{Supp}_Xz
\endcases
\tag{NP4} 
$$
and
$$
\cases
D^*(f,\lambda)/D(f,\lambda)>1\\
\Leftrightarrow \lambda=c_i(f)\not\in\Bbb Z
\text{ for some }i\in\{1,\dots,h(c(f))\}\\
\Leftrightarrow\widehat E(f,\lambda,z,Y)
\text{ has more than one root in }k\\
\qquad\text{for any $z=z(X)\in k((X))$ with $f(X^n,z(X))=0$}
\endcases
\tag{NP5} 
$$ 
and
$$
\cases
\text{$S(f,\lambda)\in\Bbb Q$ is such that:}\\
\text{if }(z,U,V,W)\in\epsilon(f,\lambda)
\text{ then $S(f,\lambda)nV\in\Bbb Z$.}
\endcases
\tag{NP6}
$$
With this preparation, what we have called {\bf Newton Polygon} (3) on 
page 334 of \cite{Ab} can be restated by saying that:
$$
\cases
\text{if }(z,U,V,W)\in\epsilon(f,\lambda)\\
\text{then ord}_Xf(X^V,z^{\ddagger}(X,U,V,W,Y))=S(f,\lambda)nV\\ 
\text{and inco}_Xf(X^V,z^{\ddagger}(X,U,V,W,Y))=\widehat A(f,\lambda,z)
\widehat E(f,\lambda,z,Y).
\endcases
\tag{NP7}
$$
In view of (NP3) and (NP5), the last line of (NP7) tells us that the 
noninteger members of the sequence $c_i(f)_{1\le i\le h(c(f))}$ 
are exactly those values of $\lambda$ for which 
inco$_X(f(X^V,z(X,U,V,W,Y))$ has more than one root in $k$; 
this then is the alternative characterization we spoke of.

Given any other
$f'=f'(X,Y)\in R^{\natural}$ of $Y$-degree $n'$, by Newton's Theorem
$$
f'(X^{n'},Y)=\prod_{1\le j\le n'}[Y-z'_{j}(X)]
\quad\text{ with }\quad z'_{j}(X)\in k((X)).
$$
Recall that on page 287 of \cite{Ab} the {\bf contact} cont$(f,f')$ of $f$
with $f'$ is defined by putting
$$
\text{cont}(f,f')
=\text{max}\{(1/n')\text{ord}_X[z_j(X^{n'})-z'_{j'}(X^n)]:
1\le j\le n\text{ and }1\le j'\le n'\}.
$$
We define the {\bf normalized contact} $\text{noc}(f,f')$ of $f$ with $f'$
by putting $\text{noc}(f,f')=(1/n)\text{cont}(f,f')$, i.e., equivalently, 
by putting
$$
\text{noc}(f,f')=\text{max}\{(1/(nn'))\text{ord}_X[z_j(X^{n'})-z'_{j'}(X^n)]:
1\le j\le n\text{ and }1\le j'\le n'\}.
$$
We note that if $f\ne f'$ then $\text{noc}(f,f')$ is a rational number, and if
$f=f'$ then $\text{noc}(f,f')=\infty$. 
We also note the {\bf Isosceles Triangle 
Property} which we shall tacitly use and which says that 
$$
\cases
f''\in R^{\natural}\Rightarrow
\text{noc}(f,f'')\ge\text{min}(\text{noc}(f,f'),\text{noc}(f',f''))\\
\text{and}\\
f''\in R^{\natural}\text{ with }\text{noc}(f,f')\ne \text{noc}(f',f'')
\Rightarrow
\text{noc}(f,f'')=\text{min}(\text{noc}(f,f'),\text{noc}(f',f'')).
\endcases
\tag{ITP}
$$

In view of the Confluence Lemmas given on pages 338-344 of \cite{Ab}
we see that
$$
\cases
\text{if }\lambda\le\lambda'=\text{noc}(f,f')\\
\text{then $p(f',\lambda)=p(f,\lambda)$, $D(f',\lambda)=D(f,\lambda)$,
$S(f',\lambda)=S(f,\lambda)$,}\\
\text{$m(f',\lambda)=m(f,\lambda)$, $c(f',\lambda)=c(f,\lambda)$, 
$t(f',\lambda)=t(f,\lambda)$,}\\
\text{and $A(f',\lambda,z')=A(f,\lambda,z)$}\\
\text{where we have chosen }z=z(X)\text{ and }z'=z'(X)\text{ in }k((X))
\text{ such that}\\
f(X^n,z(X))=0=f'(X^{n'},z'(X))\text{ and }
(1/(nn'))\text{ord}_X[z(X^{n'})-z'(X^n)]=\lambda'
\endcases
\tag{GNP1}
$$
and
$$
\cases
\text{if }\lambda<\lambda'=\text{noc}(f,f')\\
\text{then $p^*(f',\lambda)=p^*(f,\lambda)$, 
$D^*(f',\lambda)=D^*(f,\lambda)$,}\\
\text{$m^*(f',\lambda)=m^*(f,\lambda)$, $c^*(f',\lambda)=c^*(f,\lambda)$, 
$t^*(f',\lambda)=t^*(f,\lambda)$,}\\
\text{and $E(f',\lambda,z',Y)=E(f,\lambda,z,Y)$}\\
\text{where we have chosen }z=z(X)\text{ and }z'=z'(X)\text{ in }k((X))
\text{ such that}\\
f(X^n,z(X))=0=f'(X^{n'},z'(X))\text{ and }
(1/(nn'))\text{ord}_X[z(X^{n'})-z'(X^n)]=\lambda'
\endcases
\tag{GNP2}
$$
and
$$
\cases
\text{if }\lambda=\lambda'=\text{noc}(f,f')\\
\text{then }E(f',\lambda,z',Y)\text{ and }E(f,\lambda,z,Y)
\text{ do not have a common root}\\
\text{where we have chosen }z=z(X)\text{ and }z'=z'(X)\text{ in }k((X))
\text{ such that}\\
f(X^n,z(X))=0=f'(X^{n'},z'(X))\text{ and }
(1/(nn'))\text{ord}_X[z(X^{n'})-z'(X^n)]=\lambda'
\endcases
\tag{GNP3}
$$
and
$$
\cases
\text{if }\lambda\le\lambda'=\text{noc}(f,f')\text{ and }
(z,U,V,W)\in\epsilon(f,\lambda)\\
\text{then }S(f,\lambda)n'V\in\Bbb Z.
\endcases
\tag{GNP4}
$$
What we have called {\bf Generalized Newton Polygon} (6) on page 
346-347 of \cite{Ab} can now be restated by saying that:
$$
\cases
\text{if }\lambda\le\lambda'=\text{noc}(f,f')\text{ and }
(z,U,V,W)\in\epsilon(f,\lambda)\\
\text{then ord}_Xf'(X^V,z^{\ddagger}(X,U,V,W))=S(f,\lambda)n'V\\
\text{and }\text{inco}_Xf'(X^V,z^{\ddagger}(U,V,W,Y))
=\widehat A(f',\lambda,z')\widehat E(f',\lambda,z',Y)\\
\text{where we have chosen }z'=z'(X)\in k((X))\text{ such that}\\
f'(X^{n'},z'(X))=0\text{ and }
(1/(nn'))\text{ord}_X[z(X^{n'})-z'(X^n)]=\lambda'
\endcases
\tag{GNP5}
$$
and
$$
\cases
\text{if }
\lambda>\lambda'=\text{noc}(f,f')\text{ and }
(z,U,V,W)\in\epsilon(f,\lambda)\\
\text{then ord}_Xf'(X^V,z^{\ddagger}(X,U,V,W))=S(f,\lambda')n'V\\
\text{and }0\ne\text{inco}_Xf'(X^V,z^{\ddagger}(U,V,W,Y))
=\widehat A(f',\lambda',z')\widehat E(f',\lambda',z',z[n\lambda])\in k\\
\text{where we have chosen }z'=z'(X)\in k((X))\text{ such that}\\
f'(X^{n'},z'(X))=0\text{ and }
(1/(nn'))\text{ord}_X[z(X^{n'})-z'(X^n)]=\lambda'
\endcases
\tag{GNP6}
$$
and
$$
\cases
\text{if }
\lambda=\lambda'=\text{noc}(f,f')\text{ and }
f(X^n,z(X))=0\text{ with }z(X)\in k((X))\\
\text{then ord}_Xf'(X^n,z(X))=S(f,\lambda)n'n.
\endcases
\tag{GNP7}
$$
Finally we note that, for the truncations
$t(f,\lambda)$ and $t^*(f,\lambda)$, we obviously have
$$
\cases
\text{noc$(f,t(f,\lambda))\ge\lambda$ and}\\
\text{noc$(f,t^*(f,\lambda))>\lambda$.}
\endcases
\tag{GNP8}
$$

\centerline{}

\centerline{{\bf Section 4: Truncations and Buds}}

To continue discussing truncations, we let
$R^{\flat}$ be the set of 
all {\bf buds} in $R$, where by a bud we mean a pair 
$B=(\sigma(B),\lambda(B))$ with 
$\emptyset\ne\sigma(B)\subset R^{\natural}$ and 
$\lambda(B)\in\Bbb Q$ such that 
noc$(f,f')\ge\lambda(B)$ for all $f$ and $f'$ in $\sigma(B)$; 
we call $\sigma(B)$ the {\bf stem} of $B$, 
and $\lambda(B)$ the {\bf level} of $B$; we also let 
$\tau(B)=\{f\in R^{\natural}:
\text{noc}(f,f')\ge\lambda(B)\text{ for all }f'\in\sigma(B)\}$,
and we call $\tau(B)$ the {\bf flower} of $B$.
\footnote{Basically, the stem $\sigma(B)$ of a bud $B=(\sigma(B),\lambda(B))$
is a nonempty set of irreducible meromorphic curves $f(X,Y)=0$ whose fractional
meromorphic roots mutually coincide up to $X^{\lambda(B)}$, and its flower 
$\tau(B)$ is the set of all irreducible meromorphic curves whose fractional
meromorphic roots coincide with the fractional meromorphic roots of members of 
$\sigma(B)$ up to $X^{\lambda(B)}$.}
For any $f\in R^{\natural}$ and $B\in R^{\flat}$, we let
noc$(f,B)$ be the rational number defined by saying that if
$f\in\tau(B)$ then noc$(f,B)=\lambda(B)$, whereas if $f\not\in\tau(B)$
then noc$(f,B)$ equals the common value (see (ITP)) of noc$(f,f')$
as $f'$ varies in $\tau(B)$; we call noc$(f,B)$ the {\bf normalized
contact} of $f$ with $B$, and we note that: $\text{noc}(f,B)\ne\lambda(B)
\Leftrightarrow\text{noc}(f,B)<\lambda(B)\Leftrightarrow f\not\in\tau(B)$.
\footnote{Equivalently, noc$(f,B)$ can be defined by saying that,
for any $f\in R^{\natural}$ and $B\in R^{\flat}$, we have
noc$(f,B)=\text{min}\{\text{noc}(f,f'):f'\in\tau(B)\}$.}
For any $f\in R^{\natural}$ and $\lambda\in\Bbb Q$,
we let $\overline R(f,\lambda)=\{B'\in R^{\flat}:f\in\tau(B')
\text{ and }\lambda(B')=\lambda\}$; members of 
$\overline R(f,\lambda)$ may be called {\bf $\lambda$-buddies} of $f$. 
For any $f\in R^{\natural}$ and $B\in R^{\flat}$, we let
$\overline R(f,B)=\overline R(f,\text{noc}(f,B))$; members of 
$\overline R(f,B)$ may be called {\bf $B$-buddies} of $f$.
For any $B\in R^{\flat}$, we let $\overline R(B)=\{B'\in R^{\flat}:
\tau(B')=\tau(B)\text{ and }\lambda(B')=\lambda(B)\}$; members of
$\overline R(B)$ may be called {\bf buddies} of $B$.

{\bf Given any bud $B$,} by (GNP1) we see that there is a unique
nonnegative integer $p(B)$, a unique positive integer $D(B)$, a unique
rational number $S(B)$, a unique sequence of integers
$m(B)=m_i(B)_{0\le i\le p(B)+1}$ with the exception that 
$m_{p(B)+1}=\infty$, a unique sequence of rational numbers
$c(B)=c_i(B)_{1\le i\le p(B)}$, a unique member $t(B)$ of $R^{\natural}$,
a unique nonzero element $A(B)$ of $k$, and a unique nonempty set 
$\epsilon(B)$ of triples $(z,V,W)$ with $z=z(X)\in k((X))$ and 
$0<V\in\Bbb Z$ and $W\in\Bbb Z$, having the {\bf Bud Properties}
which say that
$$
\cases
\text{for every $f\in\tau(B)$, upon letting deg$_Yf=n$, we have:}\\
\text{$p(f,\lambda(B))=p(B)$, $D(f,\lambda(B))=D(B)$,
$S(f,\lambda(B))=S(B)$,}\\
\text{$m(f,\lambda(B))=m(B)$, $c(f,\lambda(B))=c(B)$, 
$t(f,\lambda(B))=t(B)$,}\\
\text{$A(f,\lambda(B),\widetilde z)=A(B)$ for all 
$\widetilde z=\widetilde z(X)\in k((X))$ with $f(X^n,\widetilde z(X))=0$,}\\
\text{and $(z,n,V,W)\mapsto(\widehat z,V,W)$ gives a surjection of
$\epsilon(f,\lambda(B))$ onto $\epsilon(B)$}\\
\text{where $\widehat z(X)=z^{\dagger}(X,n,V,W)$.}
\endcases
\tag{BP1}
$$
We call $p(B)$, $D(B)$, $S(B)$, $m(B)$, $c(B)$, $t(B)$, $A(B)$, and
$\epsilon(B)$, the {\bf position}, the {\bf degree}, the {\bf strength}, 
the {\bf characteristic sequence}, the {\bf normalized characteristic 
sequence}, the {\bf normalized truncation}, the {\bf reduced constant},
and the {\bf edge} of $B$ respectively, and we note that then for $t(B)$
we have
$$
\cases
\text{$t(B)\in\tau(B)$, deg$_Yt(B)=D(B)$, $m(t(B))=m(B)$,}\\
\text{$c(t(B))=c(B)$, $h(m(t(B)))=h(c(t(B)))=p(B)$, and}\\
\text{$\epsilon(B)=\{(z,V,W):0<V\in D(B)\Bbb Z$ and $W=\lambda(B)V\in\Bbb Z$ 
and}\\
\qquad\qquad\qquad\qquad\quad
\text{$z=z(X)\in k((X))$ with $t(B)(X^V,z(X))=0$\}.}
\endcases
\tag{BP2}
$$
{\bf Given any bud $B$,} by (BP1) and (BP2) we see that
$$
\cases
\text{for any $B'\in R^{\flat}$ we have:}\\
\text{$B'\in \overline R(B)\Leftrightarrow\overline R(B')=\overline R(B)$}\\
\qquad\qquad\quad\text{$\Leftrightarrow\tau(B')\cap\tau(B)\ne\emptyset$ and 
$\lambda(B')=\lambda(B)$}\\
\qquad\qquad\quad
\text{$\Rightarrow p(B')=p(B)$, $D(B')=D(B)$, $S(B')=S(B)$, $m(B')=m(B)$,}\\
\qquad\qquad\qquad
\text{$c(B')=c(B)$, $t(B')=t(B)$, $A(B')=A(B)$, 
and $\epsilon(B')=\epsilon(B)$}
\endcases
\tag{BP3}
$$
and by (NP1), (GNP4), (GNP5) and (GNP6) we see that
$$
\cases
\text{for any $f\in\tau(B)$, with deg$_Yf=n$, we have $n/D(B)\in\Bbb Z$,}\\
\text{and for any $(z,V,W)\in\epsilon(B)$ we have}\\
\text{ord}_Xf(X^V,z(X)+X^WY)=S(B)nV\in\Bbb Z\\
\text{and deg}_Y\text{inco}_Xf(X^V,z(X)+X^WY)=n/D(B)
\endcases
\tag{BP4}
$$
and
$$
\cases
\text{for any $f'\in R^{\natural}\setminus\tau(B)$ and 
$(z,V,W)\in\epsilon(B)$ we have}\\
\text{$0\ne\text{inco}_Xf'(X^V,z(X)+X^WY)\in k$,}\\
\text{and for any $B'\in\overline R(f',B)$,
upon letting deg$_Yf'=n'$, we have}\\
\text{ord}_Xf'(X^V,z(X)+X^WY)=S(B')n'V\in\Bbb Z.\\
\endcases
\tag{BP5}
$$

Next we let $R^{\flat *}$ be the set of all {\bf strict buds} in $R$,
where by a strict bud we mean a bud $B$ such that noc$(f,f')>\lambda(B)$
for all $f$ and $f'$ in $\sigma(B)$; we also let 
$\tau^*(B)=\{f\in R^{\natural}:\text{noc}(f,f')>\lambda(B)
\text{ for all }f'\in\sigma(B)\}$, and we
call $\tau^*(B)$ the {\bf strict flower} of $B$.
\footnote{Again, basically, the stem $\sigma(B)$ of a strict bud 
$B=(\sigma(B),\lambda(B))$ is a nonempty set of irreducible meromorphic curves 
$f(X,Y)=0$ whose fractional meromorphic roots mutually coincide thru 
$X^{\lambda(B)}$, and its strict flower $\tau^*(B)$ is the set of all 
irreducible meromorphic curves whose fractional meromorphic roots coincide 
with the fractional meromorphic roots of members of $\sigma(B)$ thru 
$X^{\lambda(B)}$.}
For any $f\in R^{\natural}$ and $\lambda\in\Bbb Q$
we let $\overline R^*(f,\lambda)=\{B'\in R^{\flat *}:f\in\tau^*(B')
\text{ and }\lambda(B')=\lambda\}$; members of $\overline R^*(f,\lambda)$ 
may be called {\bf strict $\lambda$-buddies} of $f$. 
For any $f\in R^{\natural}$ and $B\in R^{\flat}$, we let
$\overline R^*(f,B)=\overline R^*(f,\text{noc}(f,B))$; members of 
$\overline R^*(f,B)$ may be called {\bf strict $B$-buddies} of $f$.
For any $B\in R^{\flat}$, we let 
$\overline R^*(B)=R^{\flat *}\cap\overline R(B)$; members
of $\overline R^*(B)$ may be called {\bf strict buddies} of $B$. Finally,
for any $B\in R^{\flat *}$, we let 
$\overline R^{**}(B)=\{B'\in\overline R^*(B):\tau^*(B')=\tau^*(B)\}$; members 
of $\overline R^{**}(B)$ may be called {\bf doubly strict buddies} of $B$.

{\bf Given any strict bud $B$,} by (GNP2) we see that there is 
a unique nonnegative integer $p^*(B)$, a unique positive integer $D^*(B)$, 
a unique sequence of integers $m^*(B)=m^*_i(B)_{0\le i\le p^*(B)+1}$ 
with the exception that $m^*_{p^*(B)+1}=\infty$, and a unique sequence of 
rational numbers $c^*(B)=c^*_i(B)_{1\le i\le p^*(B)}$, 
a unique member $t^*(B)$ of $R^{\natural}$, a unique monic polynomial 
$E(B,Y)$ in $k[Y]$, a unique element $E_0(B)$ in $k$, 
and a unique nonempty set $\epsilon^*(B)$ of triples 
$(z,V,W)$ with $z=z(X)\in k((X))$ and $0<V\in\Bbb Z$ and $W\in\Bbb Z$, 
having the {\bf Strict Bud Properties} which say that
$$
\cases
\text{for every $f\in\tau^*(B)$, upon letting deg$_Yf=n$, we have:}\\
\text{$p^*(f,\lambda(B))=p^*(B)\ge p(B)$, 
$D^*(f,\lambda(B))=D^*(B)\in D(B)\Bbb Z$, }\\
\text{$m^*(f,\lambda(B))=m^*(B)$, $c^*(f,\lambda(B))=c^*(B)$, 
$t^*(f,\lambda(B))=t^*(B)$,}\\
\text{$E(f,\lambda(B),\widetilde z,Y)=E(B,Y)=Y^{D^*(B)/D(B)}-E_0(B)$}\\
\text{for all $\widetilde z=\widetilde z(X)\in k((X))$ 
with $f(X^n,\widetilde z(X))=0$,}\\
\text{and $(z,n,V,W)\mapsto(\widehat z,V,W)$ gives a surjection of
$\epsilon(f,\lambda(B))$ onto $\epsilon^*(B)$}\\
\text{where $\widehat z(X)=z^{\dagger *}(X,n,V,W)$.}
\endcases
\tag{SBP1}
$$
We call $p^*(B)$, $D^*(B)$, $m^*(B)$, $c^*(B)$, 
$t^*(B)$, $E(B,Y)$, $E_0(B)$, and $\epsilon^*(B)$, 
the {\bf strict position}, the {\bf strict degree}, 
the {\bf strict characteristic sequence}, the 
{\bf strict normalized characteristic sequence},
the {\bf strict normalized truncation}, 
the {\bf reduced polynomial}, the {\bf polynomial constant},
and the {\bf strict edge} of $B$ respectively, and we note that then
for $t^*(B)$ we have
$$
\cases
\text{$t^*(B)\in\tau^*(B)$, deg$_Yt^*(B)=D^*(B)$, 
$m(t^*(B))=m^*(B)$,}\\
\text{$c(t^*(B))=c^*(B)$, $h(m(t^*(B)))=h(c(t^*(B)))=p^*(B)$, and}\\
\text{$\epsilon^*(B)=\{(z,V,W):0<V\in D^*(B)\Bbb Z$ and 
$W=\lambda(B)V\in\Bbb Z$ and}\\
\qquad\qquad\qquad\qquad\quad
\text{$z=z(X)\in k((X))$ with $t^*(B)(X^V,z(X))=0$\}.}
\endcases
\tag{SBP2}
$$
{\bf Given any strict bud $B$,} by (SBP1) and (SBP2) we see that
$$
\cases
\text{for any $B'\in R^{\flat *}$ we have:}\\
B'\in\overline R^{**}(B)\Leftrightarrow\overline R^{**}(B')
=\overline R^{**}(B)\\
\qquad\qquad\quad\;\;
\Leftrightarrow\tau^*(B')\cap\tau^*(B)\ne\emptyset\text{ and }  
\lambda(B')=\lambda(B)\\
\qquad\qquad\quad\;\;
\Rightarrow\text{$p^*(B')=p^*(B)$, $D^*(B')=D^*(B)$, $m^*(B')=m^*(B)$,
$c^*(B')=c^*(B)$,}\\
\qquad\qquad\qquad\;\;
\text{$t^*(B')=t^*(B)$, $E(B',Y)=E(B,Y)$, and $\epsilon^*(B')=\epsilon^*(B)$}
\endcases
\tag{SBP3}
$$
and by (NP1) and (GNP5) we see that
$$
\cases
\text{for any $f\in\tau^*(B)$, with deg$_Yf=n$,
we have $n/D^*(B)\in\Bbb Z$,}\\
\text{and for any $(z,V,W)\in\epsilon(B)$ we have}\\
\text{inco}_Xf(X^V,z(X)+X^WY)=A(B)^{n/D(B)}E(B,Y)^{n/D^*(B)}.
\endcases
\tag{SBP4}
$$

{\bf Given any bud $B$,} by (NP2), (NP4) and (GNP3) we get the 
{\bf Mixed Bud Properties} which say that
$$
\cases
\text{for any $B'\in\overline R^*(B)$
and $B''\in\overline R^*(B)$ we have:}\\
E(B',Y)\ne E(B'',Y)\Leftrightarrow\tau^*(B')\ne\tau^*(B'')\\
\qquad\qquad\quad\qquad\qquad\;\;\text{$\Leftrightarrow 
\tau^*(B')\cap\tau^*(B'')=\emptyset$}\\
\qquad\qquad\quad\qquad\qquad\;\;\text{$\Leftrightarrow 
E(B',Y)$ and $E(B'',Y)$ have no common root in $k$}
\endcases
\tag{MBP1}
$$
and
$$
\cases
\text{for any $B'\in\overline R^*(B)$ we have:}\\
E_0(B')=0\Leftrightarrow B'\in\overline R^*(t(B),B)\\
\qquad\qquad\quad\Rightarrow
\text{$p^*(B')=p(B')$, $D^*(B')=D(B)$, $m^*(B')=m(B')$,}\\
\qquad\qquad\qquad
\text{$c^*(B')=c(B')$, $t^*(B')=t(B')$, and
$\epsilon^*(B')=\epsilon(B')$}
\endcases
\tag{MBP2}
$$
and
$$
\cases
\text{for any $B'\in\overline R^*(B)\setminus\overline R^*(t(B),B)$
and $B''\in\overline R^*(B)\setminus\overline R^*(t(B),B)$ we have:}\\
\text{$p^*(B')=p^*(B'')$, $D^*(B')=D^*(B'')$, $m^*(B')=m^*(B'')$, 
and $c^*(B')=c^*(B'')$.}\\
\endcases
\tag{MBP3}
$$

\centerline{}

\centerline{{\bf Section 5: Contact Sets}}

Given any $F=F(X,Y)\in R$ of $Y$-degree $N$, we can write
$$
F=\prod_{0\le j\le\chi(F)}F_j
\quad\text{ where }\quad 
F_0=F_0(X)\in K((X))
$$
and
$$
F_j=F_j(X,Y)\in R^{\natural}
\quad\text{ with }\quad
\text{deg}_YF_j=N_j\text{ for }1\le j\le \chi(F)
$$
and $\chi(F)$ is a nonnegative integer such that: 
$\chi(F)=0\Leftrightarrow F\in k((X))$. 
\footnote{In other words, if $F\in k((X))$ then $\chi(F)=0$, whereas
if $F\not\in k((X))$ then $\chi(F)$ equals the number of irreducible
factors of $F$ in $R$.}
We define the {\bf contact set} 
$C(F)$ of $F$ by putting
$$
\align
C(F)=&\{c_i(F_j):1\le j\le\chi(F)\text{ and }
1\le i\le h(c(F_j))\text{ and }c_i(F_j)\not\in\Bbb Z\}\\
&\cup\{\text{noc}(F_j,F_{j'}):1\le j<j'\le\chi(F)\text{ with }F_j\ne F_{j'}\}.
\endalign
$$
Upon letting 
$$
N^{\natural}=\prod_{1\le j\le\chi(F)}N_j
$$
(with the usual {\bf convention} that the product of an empty family is $1$),
by Newton's Theorem we have
$$
F(X^{N^{\natural}},Y)
=F_0(X^{N^{\natural}})\prod_{1\le j\le N}[Y-z^{\natural}_j(X)]
\quad\text{ with }\quad z^{\natural}_j(X)\in k((X))
$$
and by the material on page 300 of \cite{Ab}, as an alternative
characterization of $C(F)$, we get 
$$
C(F)=\{(1/N^{\natural})\text{ord}_X[z^{\natural}_j(X)-z^{\natural}_{j'}(X)]:
1\le j<j'\le N\text{ with }z^{\natural}_j(X)\ne z^{\natural}_{j'}(X)\}.
$$
Note that
$$
C(F)=\emptyset\Leftrightarrow N_j=1\text{ for $1\le j\le\chi(F)$ and
$F_j=F_{j'}$ for }1\le j<j'\le\chi(F).
$$

Given any $G=G(X,Y)\in R$ of $Y$-degree $M$, we can write
$$
G=\prod_{0\le j\le\chi(G)}G_j
\quad\text{ where }\quad 
G_0=G_0(X)\in K((X))
$$
and
$$
G_j=G_j(X,Y)\in R^{\natural}
\quad\text{ with }\quad
\text{deg}_YG_j=M_j\text{ for }
1\le j\le \chi(G)
$$
and $\chi(G)$ is a nonnegative integer such that: 
$\chi(G)=0\Leftrightarrow G\in k((X))$. Note that now
$$
C(FG)=C(F)\cup C(G)\cup\{\text{noc}(F_j,G_{j'}):1\le j\le\chi(F)\text{ and }
1\le j'\le\chi(G)\text{ with }F_j\ne G_{j'}\}.
$$

Let 
$$
J(F,G)=J_{(X,Y)}(F,G)
$$ 
be the {\bf jacobian} of $F=F(X,Y)$ and
$G=G(X,Y)$ with respect to $X$ and $Y$, i.e., let
$$
J(F,G)=F_XG_Y-G_XF_Y
$$ 
where subscripts denote {\bf partial derivatives}.
Our aim is to produce a factorization of $J(F,G)$ in terms of the
contact set $C(FG)$.

In case of $F\ne 0$ we can write
$$
F=X^{N^{\sharp}}P+(\text{terms of }X\text{-degree }>N^{\sharp})
$$
where
$$
N^{\sharp}=\text{ord}_X(F)
\quad\text{ and }\quad
0\ne P=P(Y)=\text{inco}_X(F)\in k[Y]
\quad\text{ with }\quad
\text{deg}_Y(P)=\nu.
$$
Likewise, in case of $G\ne 0$ we can write
$$
G=X^{M^{\sharp}}Q+(\text{terms of }X\text{-degree }>M^{\sharp})
$$
where
$$
M^{\sharp}=\text{ord}_X(G)
\quad\text{ and }\quad
0\ne Q=Q(Y)=\text{inco}_X(G)\in k[Y]
\quad\text{ with }\quad
\text{deg}_Y(Q)=\mu.
$$
Now in case of $F\ne 0\ne G$ we get
$$
\aligned
J(F,G)=&X^{N^{\sharp}+M^{\sharp}-1}(N^{\sharp}PQ_Y-M^{\sharp}P_YQ)\\
&+(\text{terms of }X\text{-degree }>N^{\sharp}+M^{\sharp}-1)
\endaligned
\tag{JE1}
$$
and hence 
$$
\text{ord}_XJ(F,G)\ge N^{\sharp}+M^{\sharp}-1
\tag{JE2}
$$
and
$$
\aligned
\text{ord}_XJ(F,G)= N^{\sharp}+M^{\sharp}-1
&\Leftrightarrow N^{\sharp}PQ_Y-M^{\sharp}P_YQ\ne 0\\
&\Rightarrow\text{inco}_XJ(F,G)=N^{\sharp}PQ_Y-M^{\sharp}P_YQ.
\endaligned
\tag{JE3}
$$
 
These {\bf Jacobian Estimates} are basic in getting a factorization of 
$J(F,G)$ out of $C(FG)$ or, more precisely, out of the ``tree'' $T(FG)$ 
which, in Section 6, we shall build from $C(FG)$. 
Moreover, as we shall explain in Section 7,
most of this set-up works in getting a factorization of any 
$H\in R$ out of any tree $T$. In Section 8 we shall apply it to the situation 
when $H=J(F,G)=F_Y$ with $G=-X$. In Section 9 we shall
consider the general case of $H=J(F,G)$.

\centerline{}

\centerline{{\bf Section 6: Trees}}

By allowing the level $\lambda(B)$ of a bud $B$ to be $-\infty$ we get the 
set $R_{\infty}^{\flat}$ of all {\bf improper buds} $B$; note that any
nonempty subset of $R^{\natural}$ can be the stem $\sigma(B)$ of an improper
bud $B$; moreover, for any improper bud $B$ we have $\lambda(B)=-\infty$
and $\tau(B)=R^{\natural}$. 
We put $\widehat R^{\flat}=R^{\flat}\cup R_{\infty}^{\flat}$, and we
call a member of $\widehat R^{\flat}$ a {\bf generalized bud}.
For any $B\in\widehat R^{\flat}$ we let
$\tau^*(B)=\{f\in\tau(B):\text{noc}(f,f')>\lambda(B)\text{ for some }
f'\in\sigma(B)\}$, and we call $\tau^*(B)$ the {\bf strict flower} of $B$;
note that for any $B\in R^{\flat *}$ this definition coincides with the 
definition made earlier; also note that for any $B\in R^{\flat}_{\infty}$
we have $\tau^*(B)=R^{\natural}$. 
For any $B\in\widehat R^{\flat}$ we let
$\tau'(B)=\tau(B)\setminus\tau^*(B)$, and we call $\tau'(B)$ the 
{\bf primitive flower} of $B$.
Previously we have defined the {\bf normalized contact} noc$(f,B)$ for all
$f\in R^{\natural}$ and $B\in R^{\flat}$; now we extend this by putting 
noc$(f,B)=-\infty$ for all $f\in R^{\natural}$ and $B\in R^{\flat}_{\infty}$.
For any $f\in R^{\natural}$ and $B\in\widehat R^{\flat}$ we define
$R^*(f,B)$ to be the unique member of $\widehat R^{\flat}$ whose stem is
$\{f\}$ and whose level is noc$(f,B)$, and we call $R^*(f,B)$ the
{\bf strict $B$-friend} of $f$; note that then $R^*(f,B)$ belongs to
$R^{\flat *}$ or $R^{\flat}_{\infty}$ according as $B\in R^{\flat}$ or
$B\in R^{\flat}_{\infty}$. 
For any $B\in\widehat R^{\flat}$ we define $R^*(B)$ to be the set of all
$B'\in R^{\flat *}\cup R^{\flat}_{\infty}$ such that 
$\lambda(B')=\lambda(B)$ and $\sigma(B')=\tau^*(B')\cap\sigma(B)$,
and we call members $R^*(B)$ {\bf strict friends} of $B$; note that
$\sigma(B)=\coprod_{B'\in R^*(B)}\sigma(B')$ gives a partition of
$\sigma(B)$ into pairwise disjoint nonempty subsets.

Now the set $\widehat R^{\flat}$
is prepartially ordered by defining $B'\ge B$ 
to mean $\lambda(B')\ge \lambda(B)$ and $\tau(B')\subset\tau(B)$. 
\footnote{A set is prepartially ordered by $\ge$ means: 
$a\ge b$ and $b\ge c$ implies $a\ge c$. It is partially ordered if also:
$a\ge b$ and $b\ge a$ implies $a=b$.}
For any $B'$ and $B$ in $\widehat R^{\flat}$ we write $B'>>B$ or
$B<<B'$ to mean $B'>B$ and $\lambda(B')>\lambda(B)$, i.e., to mean
$\lambda(B')>\lambda(B)$ and $\tau(B')\subset\tau(B)$.
For any $B'>>B$ in $\widehat R^{\flat}$ we define
$R^*(B',B)$ to be the unique member of $\widehat R^{\flat}$ whose stem is
$\sigma(B')$ and whose level is $\lambda(B)$, and we call $R^*(B',B)$ the
{\bf strict $B$-friend} of $B'$; note that then $R^*(B',B)$ belongs to
$R^{\flat *}$ or $R^{\flat}_{\infty}$
according as $\lambda(B)\ne-\infty$ or $\lambda(B)=-\infty$.
For any $B'>>B$ in $\widehat R^{\flat}$ we also put
$\tau^*(B',B)=\tau^*(R^*(B',B))\setminus\tau(B')$, and we call
$\tau^*(B',B)$ the {\bf strict $B$-flower} of $B'$.

Let $\widehat R^{\sharp}$ be the set of all {\bf trees} in $R$, where by a 
tree we mean a subset $T$ of $\widehat R^{\flat}$ such that $T$ contains an 
improper bud, and for any $B'\ne B$ in $T$ with
$\lambda(B')=\lambda(B)$ we have $\tau(B')\cap\tau(B)=\emptyset$; note that 
then the prepartial order $\ge$ induces a partial order on $T$, and hence in
particular $T$ has a unique improper bud; we call this improper bud 
the {\bf root} of $T$ and denote it by $R_{\infty}(T)$; also note that
for any $B'$ and $B$ in $T$ we have: $B'>B\Leftrightarrow B'>>B$. 
For any tree $T$, we put $\Lambda(T)=\{\lambda(B):B\in T\}$ and we 
call $\Lambda(T)$ the {\bf level set} of $T$; we define the {\bf height}
$h(T)$ of $T$ by putting $h(T)=\infty$ if $\Lambda(T)$ is infinite,
and $h(T)=$ the cardinality of $\Lambda(T)$ minus $1$ if $\Lambda(T)$
is finite; moreover, in case $h(T)$ is a nonnegative integer, i.e.,
in case $\Lambda(T)$ is a finite set, we let
$l(T)=l_i(T)_{0\le i\le h(T)}$ be the strictly increasing sequence
$l_0(T)<\dots<l_{h(T)}(T)$ such that $\{l_0(T),\dots,l_{h(T)}(T)\}
=\Lambda(T)$, and we call $l(T)$ the {\bf level sequence} of $T$.
Note that a tree $T$ is {\bf finite} iff its level set $\Lambda(T)$ is finite 
and $T$ has at most a finite number of generalized buds of any given level.
We put
$$
R^{\sharp}=\text{the set of all {\bf finite trees} in }R.
$$

For any generalized bud $B$ in any tree $T$, we put 
$\pi(T,B)=\{B'\in T:B'>B\}$ and we 
call $\pi(T,B)$ the $B$-{\bf preroof} of $T$, and we put
$\rho(T,B)=\{B'\in\pi(T,B):
\text{there is no $B''\in\pi(T,B)$ with }B'>B''\}$ 
and  we call $\rho(T,B)$ the $B$-{\bf roof} of $T$. 
For any generalized bud $B$ in any tree $T$, we also put
$\tau(T,B)=\tau(B)\setminus\cup\{\tau(B'):B'\in\rho(T,B)\}$ and we call
$\tau(T,B)$ the {\bf $B$-flower} of $T$, and we put
$\tau^*(T,B)=\tau^*(B)\setminus\cup\{\tau(B'):B'\in\rho(T,B)\}$ and we call
$\tau^*(T,B)$ the {\bf strict $T$-flower} of $B$; note that then
$\tau(T,B)=\tau(B)\setminus\cup\{\tau(B'):B'\in\pi(T,B)\}$ and
$\tau^*(T,B)=\tau^*(B)\setminus\cup\{\tau(B'):B'\in\pi(T,B)\}
=\tau(T,B)\setminus\tau'(B)$.
\footnote{For printing convenience we may write
$\cup\{\tau(B'):B'\in\rho(T,B)\}$ instead of
$\cup_{B'\in\rho(T,B)}\tau(B')$, with similar notation for
$\cap$, $\sum$ and $\prod$.}

A tree $T$ is said to be {\bf strict} if for every $\lambda\in\Lambda(T)$
we have $\sigma(R_{\infty}(T))
=\cup_{B\in T^{(\lambda)}}\sigma(B)$ where $T^{(\lambda)}$ is the set
of all $B\in T$ with $\lambda(B)=\lambda$. 
Given any $\lambda\in{\Bbb Q}\cup\{-\infty\}$, by (ITP) we see that
$f\sim_{\lambda}f'$ gives an equivalence relation on $R^{\natural}$
where $f\sim_{\lambda}f'$ means noc$(f,f')\ge\lambda$. It follows that,
given any $\widehat\sigma\subset R^{\natural}$ and 
$\widehat\Lambda\subset\Bbb Q$,
there is a unique strict tree $\widehat T(\widehat\sigma,\widehat\Lambda)$ 
with $\Lambda(\widehat T(\widehat\sigma,\widehat\Lambda))
=\{-\infty\}\cup\widehat\Lambda$
such that $\sigma(R_{\infty}(\widehat T(\widehat\sigma,\widehat\Lambda)))
=\widehat\sigma$ or $\{Y\}$ according as $\widehat\sigma$ is nonempty or 
empty; we call $\widehat T(\widehat\sigma,\widehat\Lambda)$ the
{\bf $\widehat\Lambda$-tree} of $\widehat\sigma$; 
note that, if $\widehat\sigma$ is nonempty then, for every 
$\lambda\in\widehat\Lambda$, the stems of the buds of 
$\widehat T(\widehat\sigma,\widehat\Lambda)$ of level $\lambda$
are the equivalence classes of $\widehat\sigma$ under $\sim_{\lambda}$;
likewise, if $\widehat\sigma$ is empty then, for every 
$\lambda\in\widehat\Lambda$, the stem of the unique bud of
$\widehat T(\widehat\sigma,\widehat\Lambda)$ of level $\lambda$ is $\{Y\}$. 
We put
$$
R^{\sharp *}=\text{the set of all {\bf finite strict trees} in }R
$$
and we note that for any $B\in T\in R^{\sharp *}$ with $\lambda(B)=l_i$ 
for some $i<h(T)$ we have 
$\rho(T,B)=\{B'\in\ T:\lambda(B')=l_{i+1}\}$ and
$\sigma(B)=\coprod_{B'\in\rho(T,B)}\sigma(B')$ 
which is a partition of $\sigma(B)$ into pairwise disjoint nonempty subsets. 
For any $F\in R$, with its monic irreducible factors 
$F_1,\dots,F_{\chi(F)}$ as in the previous Section, we put 
$$
T(F)=\widehat T(\{F_1,\dots,F_{\chi(F)}\},C(F))
$$
and we call $T(F)$ the {\bf tree} of $F$, and we note that then
$T(F)\in\ R^{\sharp *}$. 

A tree $T'$ is a {\bf subtree} of a tree $T$ if for 
every $B'\in T'$ there exists some (and hence a unique) $B\in T$ such that 
$\sigma(B')\subset\sigma(B)$ and $\lambda(B')=\lambda(B)$.
Every tree is clearly a subtree of the {\bf universal tree}
$\widehat T(R^{\natural},\Bbb Q)$, which is a strict tree of
infinite height.
\footnote{This universal tree is like the ASHWATTHA TREE of the Bhagwad-Gita. 
The stem of its root contains the embryos of all the past, present and 
future creatures in nascent form. 
Its trunks travel upwards first comprised of large tribes and then of smaller 
and smaller clans. Its ``ultimate'' shoots reaching heaven are the individual 
souls eager to embrace their maker.} 

\centerline{}

\noindent
{\bf Remark (TR1).} Basically we are interested in comparing the tree
$T(FG)$ of the product of two members $F$ and $G$ of $R$ with the tree
$T(J(F,G))$ of their jacobian. In case of $G=-X$, this reduces to
comparing $T(F)$ with $T(F_Y)$.

\centerline{}

\noindent
{\bf Remark (TR2).} For the benefit of the readers (and ourselves) we shall 
now describe three examples of the tree $T(F)$ of various types of 
$F\in\widehat R^{\natural}$.

\centerline{}

\noindent
{\bf Example (TR3).} First, here is an example of 
$F\in\widehat R^{\natural}$ 
which is irreducible and has only one characteristic exponent, i.e., 
with $\chi(F)=1$ and $h(m(F))=1$. Namely, let
$$
1<n\in\Bbb Z\text{ and }0\ne e\in\Bbb Z\text{ with GCD}(n,e)=1
$$
and
$$
F=f=f(X,Y)=Y^n+\sum_{1\le i\le n}w_i(X)Y^{n-i}
$$
where $w_i(X)\in k((X))$ is such that
$$
\text{ord}_X w_i(X)>ie/n\text{ for }1\le i\le n-1
\text{ and }\text{ord}_Xw_n(X)=e
$$
and let $\kappa$ be the coefficient of $X^e$ in $w_n(X)$, i.e., let
$0\ne\kappa\in k$ be such that ord$_X(w_n(X)-\kappa X^e)>e$.
Then $f$ is irreducible in $\widehat R^{\natural}$, and 
we have the Newtonian Factorization
$$
f(X^n,Y)=\prod_{1\le j\le n}[Y-z_j(X)]
$$
where $z_j(X)\in k((X))$ is such that
$$
z_j(X)=\omega^j \kappa^* X^e+(\text{terms of degree $>e$ in $X$})
$$
where $\omega$ is a primitive $n$-th root of $1$ in $k$, and $\kappa^*$
is an $n$-th root of $-\kappa$ in $k$.

To see this, first note that
$f(X^n,X^eY)=X^{ne}g(X,Y)$ where
$$
g(X,Y)=Y^n+\sum_{1\le i\le n}v_i(X)Y^{n-i}
$$ 
and $v_i(X)=X^{-ie}w_i(X^n)\in k[[X]]$ is such that
$$
v_i(0)=0\text{ for }1\le i\le n-1\text{ and } v_n(0)=\kappa.
$$
Now $g(0,Y)=Y^n-\kappa^{*n}$, and hence we get the desired factorization by
applying Hensel's Lemma. Since GCD$(n,e)=1$, we see that
$f$ is irreducible in $\widehat R^{\natural}$.

The above factorization of $f$ yields $h(m(f))=1$ with
$$
m_0(f)=q_0(m(f))=s_0(q(m(f)))=r_0(q(m(f)))=n
\text{ and } d_1(m(f))=n
$$
and
$$
m_1(f)=q_1(m(f))=s_1(q(m(f)))/n=r_1(q(m(f)))=e
\text{ and } d_2(m(f))=1.
$$
Therefore
$$
C(F)=C(f)=\{c_1(f)\}
\text{ with }
c_1(f)=e/n.
$$
Hence $h(T(f))=1$ with 
$$
l_0(T(f))=-\infty\text{ and }l_1(T(f))=c_1(f)=e/n
$$
and upon letting 
$$
B_i\in\widehat R^{\flat}\text{ with }
\sigma(B_i)=\{f\}\text{ and }\lambda(B_i)=l_i(T(f))
\text{ for }0\le i\le 1
$$
we have
$$
T(F)=T(f)=\{B_0,B_1\}
$$
with
$$
D'(B_0)=0\text{ and }D'(B_1)=n-1.
$$

Note that for $F$ to be {\bf analytic}, i.e., for it to belong to 
the ring $k[[X]][Y]$, the condition $e>0$ is necessary and sufficient. 
However, for $F$ to be {\bf pure meromorphic}, i.e., 
for it to belong to the ring
$k[X^{-1}][Y]$, i.e., for the existence of $\Phi(X,Y)\in k[X,Y]$ with
$F(X,Y)=\Phi(X^{-1},Y)$, the condition $e<0$ is necessary but not sufficient,
As a specific illustration of the analytic case we may take
$(n,e)=(4,5)$ and $(w_1(X),\dots,w_{n-1}(X),w_n(X))=(0,\dots,0,X^5)$, 
giving us $F(X,Y)=Y^4+X^5$. Similarly, as a specific illustration of the pure 
meromorphic case we may take 
$(n,e)=(4,-3)$ and $(w_1(X),\dots,w_{n-1}(X),w_n(X))=(0,\dots,0,X^{-3})$, 
giving us $F(X,Y)=Y^4+X^{-3}$, i.e., $F(X,Y)=\Phi(X^{-1},Y)$ with
$\Phi(X,Y)=Y^4+X^3$.

\centerline{}

\noindent
{\bf Example (TR4).} Next, here is an example of 
$F\in\widehat R^{\natural}$ 
which is irreducible and has two characteristic exponents, i.e., 
with $\chi(F)=1$ and $h(m(F))=2$. Namely, let
$$
F=f=f(X,Y)=(Y^2-X^{2a+1})^2-X^{3a+b+2}Y\text{ with }
a\in\Bbb Z\text{ and }
0\le b\in\Bbb Z.
$$
Then $f$ is irreducible in $\widehat R^{\natural}$, and
we have the Newtonian Factorization
$$
f(X^4,Y)=\prod_{1\le j\le 4}[Y-z_j(X)]
$$
where $z_j(X)\in k((X))$ is such that
$$
z_j(X)=(\iota^jX)^{4a+2}+\frac{1}{2}(\iota^jX)^{4a+2b+3}
+(\text{terms of degree $>4a+2b+3$ in $X$})
$$
where $\iota$ is a primitive $4$-th root of $1$ in $k$ (e.g., $\iota=$ the
usual $i$).

To see this, first note that
$f(X^4,X^{4a+2}Y)=X^{16a+8}g(X,Y)$ where
$$
g(X,Y)=(Y^2-1)^2-X^{4b+2}Y.
$$
Now for $\jmath=1$ or $-1$,
upon letting $g_{\jmath}(X,Y)=g(X,Y+\jmath)$ we have
$$
g_{\jmath}(X,Y)=(2\jmath+Y)^2Y^2-\jmath(1+\jmath Y)X^{4b+2}
$$ 
and hence (say by the Binomial Theorem) we get
$$
g_{\jmath}(X,Y)=[(2\jmath+Y)Y-\jmath^*\theta(Y)X^{2b+1}]
[(2\jmath+Y)Y+\jmath^*\theta(Y)X^{2b+1}]
$$
where 
$$\theta(Y)=1+(\jmath Y/2)
-\sum_{i=2}^{\infty}1\times 3\times\dots\times(2i-3)\times(-\jmath Y/2)^i/i!
$$
and $\jmath^*=\iota^2$ or $\iota$ according as $\jmath=1$ or $-1$, 
and therefore (say by the Weierstrass Preparation Theorem) we have
$$
g(X,Y)=\prod_{1\le j\le 4}[Y-y_j(X)]
$$ 
where $y_j(X)\in k[[X]]$ is such that
$$
y_j(X)=(-1)^j\left[1+\frac{1}{2}(\iota^jX)^{2b+1}
+(\text{terms of degree $>2b+1$ in $X$})\right].
$$
Since
$f(X^4,X^{4a+2}Y)=X^{16a+8}g(X,Y)$, we get the above factorization of 
$f(X^4,Y)$. Since the GCD of $4$ with the support of $z_1(X)$ is $1$, we
conclude that $f$ is irreducible in $\widehat R^{\natural}$, i.e.,
$\chi(F)=1$. 

The above factorization of $f$ yields $h(m(f))=2$ with
$$
m_0(f)=q_0(m(f))=s_0(q(m(f)))=r_0(q(m(f)))=4
\text{ and }d_1(m(f))=4
$$
and
$$
m_1(f)=q_1(m(f))=s_1(q(m(f)))/4=r_1(q(m(f)))=4a+2
\text{ and }d_2(m(f))=2
$$
and
$$
m_2(f)=4a+2b+3\text{ and }q_2(m(f))=2b+1
$$
and
$$
s_2(q(m(f)))=16a+4b+10\text{ and }r_2(q(m(f)))=8a+2b+5\text{ and }
d_3(m(f))=1.
$$
Therefore
$$
C(F)=C(f)=\{c_1(f),c_2(f)\}
$$
with
$$
c_1(f)=(2a+1)/2\text{ and }c_2(f)=(4a+2b+3)/4.
$$
Hence $h(T(f))=2$ with $l_0(T(f))=-\infty$ and
$$
l_1(T(f))=c_1(f)=(2a+1)/2\text{ and }l_2(T(f))=c_2(f)=(4a+2b+3)/4
$$
and upon letting 
$$
B_i\in\widehat R^{\flat}\text{ with }
\sigma(B_i)=\{f\}\text{ and }\lambda(B_i)=l_i(T(f))
\text{ for }0\le i\le 2
$$
we have
$$
T(F)=T(f)=\{B_0,B_1,B_2\}.
$$
with $D'(B_0)=0$ and 
$$
D'(B_1)=1\text{ and }D'(B_2)=2.
$$

As a specific illustration of the analytic case we may take $(a,b)=(1,0)$,
giving us $F(X,Y)=(Y^2-X^3)^2-X^5Y$. Similarly, as a specific illustration 
of the pure meromorphic case we may take $(a,b)=(-1,1)$,
giving us $F(X,Y)=(Y^2-X^{-1})^2-Y$, i.e., $F(X,Y)=\Phi(X^{-1},Y)$ with
$\Phi(X,Y)=(Y^2-X)^2-Y$. 
Note that this $\Phi$ is a {\bf variable} in the sense that
$k[X,Y]=k[\Phi,\Psi]$ for some $\Psi$ in $k[X,Y]$; in our situation
we can take $\Psi(X,Y)=Y^2-X$.

\centerline{}

\noindent
{\bf Example (TR5).} Finally, here is an example of 
$F\in\widehat R^{\natural}$ which has two factors, i.e., with $\chi(F)=2$. 
Namely, let
$$
0<n\in\Bbb Z
\text{ and }
a\in\Bbb Z
\text{ and }
0\le b\in\Bbb Z
$$
and
$$
F=F(X,Y)=Y^{n+2}+\sum_{2\le i\le n+2}u_i(X)Y^{n+2-i}
$$
where $u_i(X)\in k((X))$ is such that
$$
\text{ord}_Xu_i(X)\ge i(a+1)\text{ for }3\le i\le n+1
$$
and
$$
\text{ord}_Xu_2(X)=2a+1\text{ and ord}_Xu_{n+2}(X)=(n+2)(a+1)+b
$$
and let $0\ne\kappa'\in k$ and $0\ne\kappa\in k$ be the coefficients of
$X^{2a+1}$ and $X^{(n+2)(a+1)+b}$ in $u_2(X)$ and $u_{n+2}(X)$
respectively.
Then
$$
F(X,Y)=f(X,Y)f'(X,Y) \text{ with }f(X,Y)\ne f'(X,Y)
$$
where
$$
f(X,Y)=Y^n+\sum_{1\le i\le n}w_i(X)Y^{n-i}\in\widehat R^{\natural}
\text{ and }w_i(X)\in k((X))
$$
with
$$
\cases
\text{ord}_Xw_i(X)>ie/n\text{ for $1\le i\le n-1$ and ord}_Xw_n(X)=e+b\\
\text{for the integer $e=na+n+1$ for which GCD}(n,e)=1\\
\endcases
$$
and
$$
f'(X,Y)=Y^2+\sum_{1\le i\le 2}w'_i(X)Y^{n-i}\in R^{\natural}
\text{ and }w'_i(X)\in k((X))
$$
with
$$
\cases
\text{ord}_Xw'_1(X)>e'/2\text{ and ord}_Xw'_2(X)=e'\\
\text{for the integer $e'=2a+1$ for which GCD}(2,e')=1\\
\endcases
$$
and $0\ne\kappa'\in k$ and $0\ne \kappa/\kappa'\in k$ are the 
coefficients of $X^{e'}$ and $X^{e+b}$ in $w'_2(X)$ and $w_n(X)$
respectively.
Moreover, if $b=0$ then we also have $f(X,Y)\in R^{\natural}$.

To see this, first note that
$F(X,X^aY)=X^{na+2a}g(X,Y)$ where
$$
g(X,Y)=Y^{n+2}+\sum_{2\le i\le n+2}v_i(X)Y^{n+2-i}
$$
and $v_i(X)=X^{-ia}u_i(X)\in k[[X]]$ is such that
$$
\text{ord}_Xv_i(X)\ge i\text{ for }3\le i\le n+1
$$
and
$$
\text{ord}_Xv_2(X)=1\text{ and ord}_Xv_{n+2}(X)=n+2+b
$$
and $0\ne\kappa'\in k$ and $0\ne\kappa\in k$ are the coefficients of
$X$ and $X^{n+2+b}$ in $v_2(X)$ and $v_{n+2}(X)$ respectively.
Now the initial form $g(X,Y)$ is $\kappa' XY^n$ which
factors into the coprime factors $\kappa' X$ and $Y^n$, and hence
by the Tangent Lemma incarnation of Hensel's Lemma 
(cf. pages 140-141 of Abhyankar's 1990 AMS book 
``Algebraic Geometry For Scientists and Engineers'') we can find $\phi(X,Y)$ and
$\phi'(X,Y)$ in $k[[X,Y]]$ such that
$g(X,Y)=\phi(X,Y)\phi'(X,Y)$ and
$$
\phi(X,Y)=Y^n+\phi_{n+1}(X,Y)+(\text{terms of degree $>n+1$ in $X$ and $Y$})
$$ 
and
$$
\phi'(X,Y)=\kappa' X+\phi'_2(X,Y)+(\text{terms of degree $>2$ in $X$ and $Y$})
$$
where $\phi_{n+1}(X,Y)\in k[X,Y]$ is homogeneous of degree $n+1$
and $\phi'_{2}(X,Y)\in k[X,Y]$ is homogeneous of degree $2$  
(with the understanding that the zero polynomial is homogeneous of any degree).
Comparing terms of degree $n+2$ in the equation 
$g(X,Y)=\phi(X,Y)\phi'(X,Y)$ we get
$$
\kappa' X\phi_{n+1}(X,Y)+Y^n\phi'_2(X,Y)
=Y^{n+2}+\sum_{2\le i\le n+2}\kappa_i X^iY^{n+2-i}
$$
where $\kappa_2\in k$ is the coefficient of $X^2$ in $v_2(X)-\kappa' X$,
and $\kappa_i\in k$ is the coefficient
of $X^i$ in $v_i(X)$ for $3\le i\le n+2$.
Successively putting $X=0$ and $Y=0$ in the above equation we see that
$\phi'_2(0,Y)=Y^2$ and $\phi_{n+1}(X,0)=\kappa_{n+2} X^{n+1}$. Therefore,
in view of the Weierstrass Preparation Theorem, we can find
$\theta(X,Y)$ and $\theta'(X,Y)$ in $k[[X,Y]]$ with
$\theta(0,0)\ne 0\ne\theta'(0,0)$ such that upon letting
$\widetilde f(X,Y)=\theta(X,Y)\phi(X,Y)$ and
$\widetilde f'(X,Y)=\theta'(X,Y)\phi'(X,Y)$ we have 
$g(X,Y)=\widetilde f(X,Y)\widetilde f'(X,Y)$
and
$$
\widetilde f(X,Y)=Y^n+\sum_{1\le i\le n}\widetilde w_i(X)Y^{n-i}
\text{ and }\widetilde w_i(X)\in k[[X]]
$$
with
$$
\text{ord}_X\widetilde w_i(X)>i(n+1)/n
\text{ for $1\le i\le n-1$ and ord}_X\widetilde w_n(X)=n+1+b
$$
and
$$
\widetilde f'(X,Y)=Y^2+\sum_{1\le i\le 2}\widetilde w'_i(X)Y^{2-i}
\text{ and }\widetilde w'_i(X)\in k[[X]]
$$
with
$$
\text{ord}_X\widetilde w'_1(X)>1/2
\text{ and ord}_X\widetilde w'_2(X)=1
$$
and $0\ne\kappa'\in k$ and $0\ne\kappa/\kappa'\in k$ are the coefficients 
of $X$ and $X^{n+1+b}$ in $\widetilde w'_2(X)$ and $\widetilde w_n(X)$
respectively.
Now upon letting $f(X,Y)=X^{na}\widetilde f(X,X^{-a}Y)$ and
$f'(X,Y)=X^{2a}\widetilde f'(X,X^{-a}Y)$, 
we get the desired factorization of $F(X,Y)$. 
Since $(n,e+b)\ne (2,e')$, we also get $f\ne f'$.
By (TR3) it follows that $f'$ is irreducible in $\widehat R^{\natural}$,
and if $b=0$ then so is $f$.

Now assuming $b=0$ and $n>1$, 
in view of (TR3), the factorization of $F$ tells us that 
$h(T(F))=2$ with $l_0(T(F))=-\infty$ and 
$$
l_1(T(F))=a+(1/2)\text{ and }l_2(T(F))=a+1+(1/n)
$$
and upon letting
$$
\cases
B_0\in\widehat R^{\flat}\text{ with }
\sigma(B_0)=\{f,f'\}\text{ and }\lambda(B_0)=l_0(T(F)),\\
\text{and }B_1\in\widehat R^{\flat}\text{ with }
\sigma(B_1)=\{f,f'\}\text{ and }\lambda(B_1)=l_1(T(F)),\\
\text{and }B_2\in\widehat R^{\flat}\text{ with }
\sigma(B_2)=\{f\}\text{ and }\lambda(B_2)=l_2(T(F)),\\
\text{and }B'_2\in\widehat R^{\flat}\text{ with }
\sigma(B'_2)=\{f'\}\text{ and }\lambda(B'_2)=l_2(T(F)),\\
\endcases
$$ 
we have
$$
T(F)=\{B_0,B_1,B_2,B'_2\}
$$
with $D'(B_0)=0$ and 
$$
D'(B_1)=2\text{ and }D'(B_2)=n-1\text{ and }D'(B'_2)=0.
$$

As a specific illustration of the pure meromorphic case, taking $a=-1$ and
$(u_2(X),u_{n+2}(X))=(\kappa' X^{-1},\kappa)$ with $\kappa'\ne 0\ne\kappa$ 
in $k$ and
$u_i(X)=\kappa_i\in k$ for $3\le i\le n+1$, we get $F(X,Y)=\Phi(X^{-1},Y)$ 
where 
$$
\Phi(X,Y)=Y^{n+2}+\kappa' XY^n+\kappa
+\sum_{3\le i\le n+1}\kappa_i Y^{n+2-i}\in k[X,Y]
$$
with
$$
0\ne\kappa'\in k\text{ and }0\ne\kappa\in k
\text{ and }\kappa_i\in k\text{ for }3\le i\le n+1.
$$

\centerline{}

\centerline{{\bf Section 7: Slices}}

Given any $H=H(X,Y)\in R$ of $Y$-degree $O$, we can write 
$$
H=\prod_{0\le j\le\chi(H)}H_j
\quad\text{ where }\quad 
H_0=H_0(X)\in k((X))
\tag{SP1}
$$
and
$$
H_j=H_j(X,Y)\in R^{\natural}\quad\text{ with }\quad
\text{deg}_YH_j=O_j\text{ for }1\le j\le \chi(H)
\tag{SP2}
$$
and $\chi(H)$ is a  nonnegative integer such that: 
$\chi(H)=0\Leftrightarrow H\in k((X))$.
Now
$$
H=H_0H_{\infty}\quad\text{ with }\quad H_{\infty}=\prod_{1\le j\le\chi(H)}H_j
\tag{SP3}
$$
where we note that $H_{\infty}\in\widehat R^{\natural}$,
and we call $H_{\infty}$ the {\bf monic part} of $H$.

We put
$$
\Omega_B(H)=\prod_{1\le j\le\chi(H)\text{ with }H_j\in\tau(B)}H_j
\quad\text{ for all }B\in\widehat R^{\flat}
\tag{SP4}
$$
and we call $\Omega_B(H)$ the $B$-{\bf slice} of $H$,  
and we note that then $\Omega_B(H)\in\widehat R^{\natural}$, 
and we recall that
$$
\cases
\text{for all }B\in\widehat R^{\flat}\text{ we have}\\
\tau(B)=\{f\in R^{\natural}:\text{noc}(f,f')\ge\lambda(B)\text{ for all }
f'\in\sigma(B)\}.
\endcases
\tag{SP5}
$$
We also put 
$$
\Omega_{(T,B)}(H)
=\prod_{1\le j\le\chi(H)\text{ with }H_j\in\tau(T,B)}H_j
\quad\text{ for all }B\in T\in R^{\sharp}
\tag{SP6}
$$
and we call $\Omega_{(T,B)}(H)$ the $(T,B)$-{\bf slice} of $H$, 
and we note that then $\Omega_{(T,B)}(H)\in\widehat R^{\natural}$,
and we recall that
$$
\cases
\text{for all }B\in T\in\widehat R^{\sharp}\text{ we have}\\
\tau(T,B)=\tau(B)\setminus\cup_{B'\in\rho(T,B)}\tau(B')
=\tau(B)\setminus\cup_{B'\in\pi(T,B)}\tau(B')\\
\text{where }\pi(T,B)=\{B'\in T:B'>B\}\\
\text{and }
\rho(T,B)=\{B'\in \pi(T,B):
\text{there is no }B''\in\pi(T,B)\text{ with }B'>B''\}.
\endcases
\tag{SP7}
$$

Clearly we have the {\bf Slice Properties} which say that
$$
H_{\infty}=\Omega_B(H)\quad\text{ for all }B\in R_{\infty}^{\flat}
\tag{SP8}
$$
and
$$
\Omega_B(H)
=\Omega_{(T,B)}(H)\prod_{B'\in\rho(T,B)}\Omega_{B'}(H)
\quad\text{ for all }B\in T\in R^{\sharp}
\tag{SP9}
$$
and hence
\footnote{In the innocent looking formula (SP10), there
is more than meets the eye. Indeed it is the central theme of the paper.
It says that any finite tree $T$ gives rise to a factorization of the monic 
part $H_{\infty}$ of any meromorphic curve $H$ into the pairwise 
coprime monic factors $\Omega_{(T,B)}(H)$ with $B$ varying in $T$. 
Formula (SP20) gives a further factorization of $\Omega_{(T,B)}$
into the two coprime monic factors $\Omega'_B(H)$ and $\Omega^*_{(T,B)}(H)$.
When the finite tree $T$ is strict, formula (SP30) gives a still further
factorization of $\Omega^*_{(T,B)}(H)$ into the pairwise coprime monic
factors $\Omega^*_{(B',B)}(H)$. Item (SP50) gives a condition
for the factorization of $H_{\infty}$ to consist only of the factors
$\Omega'_B(H)$, and Item (SP80) gives a companion to this condition. 
The remaining Items (SP1)-(SP9), (SP11)-(SP19), (SP21)-(SP29), (SP31)-(SP49), 
and (SP51)-(SP79), give us details about these factors, such as their 
$Y$-degrees, and hence in particular the information as to which of these 
factors are trivial (i.e., are reduced to $1$) and which are not. Out of these 
Items, the most noteworthy are labelled as (SP40), (SP60), and (SP70).
Now roughly speaking, $\Omega_B(H)$ collects together those irreducible 
monic factors of $H$ whose normalized contact with members of $\sigma(B)$ 
is at least $\lambda(B)$, and out of these only those are kept in 
$\Omega'_B(H)$ whose normalized contact with members of $\sigma(B)$ is 
exactly $\lambda(B')$, while the remaining are put in $\Omega^*_B(H)$.
A similar description prevails for $\Omega_{(T,B)}(H)$, 
$\Omega^*_{(T,B)}(H)$, and $\Omega^*_{(B',B)}(H)$.
As we shall see in the next two Sections, more details about
these factorizations can be given when $T$ and $H$ are somehow related.}
$$
H_{\infty}=\prod_{B\in T}\Omega_{(T,B)}(H)
\quad\text{ for all }T\in R^{\sharp}
\tag{SP10}
$$
where 
$$
\cases
\text{for all $B\in T\in R^{\sharp}$ we have}\\ 
\text{deg}_Y\Omega_{(T,B)}(H)
=\text{deg}_Y\Omega_B(H)
-\sum_{B'\in\rho(T,B)}\text{deg}_Y\Omega_{B'}(H).
\endcases
\tag{SP11}
$$
By (BP4) we also see that
$$
\cases
\text{for all }B\in R^{\flat}\text{ and }(z,V,W)\in\epsilon(B)
\text{ and }f\in\tau(B)\text{ we have}\\
\text{deg}_Yf=D(B)\text{deg}_Y\text{inco}_Xf(X^V,z(X)+X^WY)
\in D(B)\Bbb Z
\endcases
\tag{SP12}
$$
and hence by (BP5) we see that
$$
\cases
\text{for all $B\in R^{\flat}$ and $(z,V,W)\in\epsilon(B)$,
upon letting $H_B=\Omega_B(H)$,}\\
\text{we have that }0\ne\text{inco}_XH_B(X^V,z(X)+X^WY)\in k[Y]\\
\text{with deg}_Y\Omega_B(H)=D(B)\text{deg}_Y\text{inco}_XH_B(X^V,z(X)+X^WY)
\in D(B)\Bbb Z\\
\text{and inco}_XH(X^V,z(X)+X^WY)
=\mu~\text{inco}_XH_B(X^V,z(X)+X^WY)\\
\text{where $\mu\in k$ is such that: $\mu=0\Leftrightarrow H=0$.}
\endcases
\tag{SP13}
$$

The factorization (SP10) can be refined further. To see this we first put
$$
\Omega'_B(H)=\prod_{1\le j\le\chi(H)\text{ with }H_j\in\tau'(B)}H_j
\quad\text{ for all }B\in\widehat R^{\flat}
\tag{SP14}
$$
and we call $\Omega'_B(H)$ the {\bf primitive $B$-slice} of $H$,  
and we note that then $\Omega'_B(H)\in\widehat R^{\natural}$,
and we recall that
$$
\cases
\text{for all }B\in\widehat R^{\flat}\text{ we have}\\
\tau'(B)=\{f\in\tau(B):\text{noc}(f,f')=\lambda(B)\text{ for all }
f'\in\sigma(B)\}.
\endcases
\tag{SP15}
$$
Next we put
$$
\Omega^*_B(H)=\prod_{1\le j\le\chi(H)\text{ with }H_j\in\tau^*(B)}H_j
\quad\text{ for all }B\in\widehat R^{\flat}
\tag{SP16}
$$
and we call $\Omega^*_B(H)$ the {\bf strict $B$-slice} of $H$,  
and we note that then $\Omega^*_B(H)\in\widehat R^{\natural}$,
and we recall that
$$
\cases
\text{for all }B\in\widehat R^{\flat}\text{ we have}\\
\tau^*(B)=\tau(B)\setminus\tau'(B)\\
\;\;\;\qquad=\{f\in\tau(B):\text{noc}(f,f')>\lambda(B)\text{ for some }
f'\in\sigma(B)\}.
\endcases
\tag{SP17}
$$
We also put
$$
\Omega^*_{(T,B)}(H)
=\prod_{1\le j\le\chi(H)\text{ with }H_j\in\tau^*(T,B)}H_j
\quad\text{ for all }B\in T\in R^{\sharp}
\tag{SP18}
$$
and we call $\Omega^*_{(T,B)}(H)$ the {\bf strict $(T,B)$-slice} of $H$,  
and we note that then $\Omega^*_{(T,B)}(H)\in\widehat R^{\natural}$,
and we recall that
$$
\cases
\text{for all }B\in T\in\widehat R^{\sharp}\text{ we have}\\
\tau^*(T,B)=\tau(T,B)\setminus\tau'(B)\\
\qquad\qquad=\tau^*(B)\setminus\cup_{B'\in\rho(T,B)}\tau(B')
=\tau^*(B)\setminus\cup_{B'\in\pi(T,B)}\tau(B').
\endcases
\tag{SP19}
$$

Now clearly
$$
\Omega_{(T,B)}(H)
=\Omega'_B(H)\Omega^*_{(T,B)}(H)
\quad\text{ for all }B\in T\in R^{\sharp}
\tag{SP20}
$$
where
$$
\cases
\text{for all $B\in R^{\flat}$ we have}\\ 
\text{deg}_Y\Omega'_B(H)
=\text{deg}_Y\Omega_B(H)-\text{deg}_Y\Omega^*_B(H).
\endcases
\tag{SP21}
$$
and
$$
\cases
\text{for all $B\in T\in R^{\sharp}$ we have}\\ 
\text{deg}_Y\Omega^*_{(T,B)}(H)
=\text{deg}_Y\Omega^*_B(H)
-\sum_{B'\in\rho(T,B)}\text{deg}_Y\Omega_{B'}(H).
\endcases
\tag{SP22}
$$
To describe the above $Y$-degrees more precisely, given any
$z=z(X)\in k((X))$, $0<V\in\Bbb Z$, and $W\in\Bbb Z$, we define the
{\bf modified $X$-initial-coefficient} of $H$ relative to $[z,V,W]$, 
to be denoted by minco$_X[z,V,W](H)$,
by putting
$$
\text{minco}_X[z,V,W](H)=\text{inco}_X H(X^V,z(X)+X^WY)
$$
and, given any $\widehat\sigma\subset R$, we define the
{\bf strict $X$-initial-coefficient} of $(H,\widehat\sigma)$ relative to 
$[z,V,W]$, to be denoted by sinco$_X[z,V,W](H,\widehat\sigma)$, and the
{\bf primitive $X$-initial-coefficient} of $(H,\widehat\sigma)$ relative to 
$[z,V,W]$, to be denoted by pinco$_X[z,V,W](H,\widehat\sigma)$, 
by saying that
$$
\cases
\text{if $H=0$ then}\\
\text{we have sinco}_X[z,V,W](H,\widehat\sigma)=1
=\text{pinco}_X[z,V,W](H,\widehat\sigma)\\
\endcases
$$
whereas
$$
\cases
\text{if $H\ne 0$ then, upon letting}\\
\text{minco}_X[z,V,W](H)=\mu_0\prod_{1\le i\le\nu}(Y-\mu_i)
\text{ with }0\ne\mu_0\in k\text{ and }\mu_i\in k\\
\text{and }\Theta_f(Y)=\text{minco}_X[z,V,W](f)
\text{ for all }f\in\widehat\sigma\\
\text{and }\sigma^*=\{i\in\{1,\dots,\nu\}:\Theta_f(\mu_i)=0
\text{ for some }f\in\widehat\sigma\}\\ 
\text{and }\sigma'=\{i\in\{1,\dots,\nu\}:\Theta_f(\mu_i)\ne 0
\text{ for all }f\in\widehat\sigma\},\\ 
\text{we have sinco}_X[z,V,W](H,\widehat\sigma)=
\prod_{i\in\sigma^*}(Y-\mu_i)\\
\text{and pinco}_X[z,V,W](H,\widehat\sigma)=
\prod_{i\in\sigma'}(Y-\mu_i).
\endcases
$$
Recall that
$$
\cases
\text{for all $B\in R^{\flat}$ we have }0\ne A(B)\in k,\\
\text{and for all $B\in R^{\flat *}$ we have }\\
E(B,Y)=Y^{D^*(B)/D(B)}-E_0(B)\text{ with }E_0(B)\in k.
\endcases
\tag{SP23}
$$
Now by (MBP2) we see that
$$
\cases
\text{for any }B\in R^{\flat *}\text{ we have:}\\ 
E_0(B)=0\Leftrightarrow B\in R^*(t(B),B)\\
\quad\qquad\qquad\Rightarrow D^*(B)=D(B)\text{ and }t^*(B)=t(B)\\
\qquad\qquad\qquad
\text{ and }\lambda(B)\ne c_i(B)\text{ for }1\le i\le p^*(B).
\endcases
\tag{SP24}
$$
and by (MBP1) we see that
$$
\cases
\text{for any $B'$ and $B''$ in $R^{\flat *}$ with 
$\tau(B')=\tau(B'')$ and $\lambda(B')=\lambda(B'')$ we have:}\\
\tau^*(B')\cap\tau^*(B'')=\emptyset\Leftrightarrow
E(B',Y)\text{ and }E(B'',Y)\text{ have no common root in }k
\endcases
\tag{SP25}
$$
and by (SBP4) we see that
$$
\cases
\text{for all }B\in R^{\flat *}\text{ with }(z,V,W)\in\epsilon(B)
\text{ and }f\in\tau^*(B)\text{ with deg}_Yf=n\\
\text{we have minco}_X[z,V,W](f)=A(B)E(B,Y)^{n/D^*(B)}
\endcases
\tag{SP26}
$$
By (SP25) and (SP26) we conclude that
$$
\cases
\text{for all }B\in R^{\flat}\text{ and }
(z,V,W)\in\epsilon(B)\text{ we have}\\
\text{that sinco$_X[z,V,W](H,\sigma(B))\in k[Y]$ is monic in $Y$}\\ 
\text{with deg}_Y\Omega^*_B(H)=D(B)\text{deg}_Y
\text{sinco}_X[z,V,W](H,\sigma(B))\in D(B)\Bbb Z\\
\text{and pinco$_X[z,V,W](H,\sigma(B))\in k[Y]$ is monic in $Y$}\\ 
\text{with deg}_Y\Omega'_B(H)=D(B)\text{deg}_Y
\text{pinco}_X[z,V,W](H,\sigma(B))\in D(B)\Bbb Z.\\
\endcases
\tag{SP27}
$$

The factorization (SP20) can be refined still further when the finite
tree $T$ is strict. To see this we put
$$
\Omega^*_{(B',B)}(H)
=\prod_{1\le j\le\chi(H)\text{ with }H_j\in\tau^*(B',B)}H_j
\quad\text{ for all }B'>>B\text{ in }\widehat R^{\flat}
\tag{SP28}
$$
and we call $\Omega^*_{(B',B)}(H)$ the {\bf strict $(B',B)$-slice} of $H$,  
and we note that then $\Omega^*_{(B',B)}(H)\in\widehat R^{\natural}$,
and we recall that
$$
\cases
\text{for all }B'>>B\text{ in }\widehat R^{\flat}\text{ we have}\\
\tau^*(B',B)=\tau^*(R^*(B',B))\setminus\tau(B')\text{ where}\\
R^*(B',B)\in\widehat R^{\flat}\text{ is given by }
\sigma(R^*(B',B))=\sigma(B')\text{ and }\lambda(R^*(B',B))=\lambda(B).
\endcases
\tag{SP29}
$$

Now clearly
$$
\Omega^*_{(T,B)}(H)
=\prod_{B'\in\rho(T,B)}\Omega^*_{(B',B)}(H)
\quad\text{ for all }B\in T\in R^{\sharp *}
\tag{SP30}
$$
where
$$
\cases
\text{for all $B\in T\in\widehat R^{\sharp}$ and $B'\in\rho(T,B)$ we have }
B'>>B \text{ in }\widehat R^{\flat},\\
\text{and in turn for all $B'>>B$ in $\widehat R^{\flat}$ we have}\\
\text{deg}_Y\Omega^*_{(B',B)}(H)
=\text{deg}_Y\Omega^*_{R^*(B',B)}(H)-\text{deg}_Y\Omega_{B'}(H)
\endcases
\tag{SP31}
$$
and
$$
\cases
\text{for all }B\in T\in R^{\sharp *}\text{ with }\lambda(B)=l_{h(T)}
\text{ we have }\rho(T,B)=\emptyset,\\
\text{whereas for all }B\in T\in R^{\sharp *}\text{ with }\lambda(B)=l_i
\text{ for some }i<h(T)\text{ we have}\\
\rho(T,B)=\{B'\in T:\lambda(B')=\lambda_{i+1}\}\text{ and }
\sigma(B)=\coprod_{B'\in\rho(T,B)}\sigma(B')\\
\text{which is a partition of $\sigma(B)$
into pairwise disjoint nonempty subsets.}
\endcases
\tag{SP32}
$$

To get more information about $\Omega_B(H)$, first we recall that
$$
\cases
\text{for any $f\in R^{\natural}$ and $B\in\widehat R^{\flat}$,}\\
R^*(f,B)\text{ is the unique member of }R^{\flat *}\cup R^{\flat}_{\infty}
\text{ with }\sigma(R^*(f,B))=\{f\}\\
\text{such that }
\lambda(R^*(f,B))=\text{min}\{\text{noc}(f,f'):f'\in\tau(B)\}
\endcases
\tag{SP33}
$$
and we put
$$
S(H,B)=\cases
\text{ord}_XH_0(X)+\sum_{1\le j\le\chi(H)}O_jS(R^*(H_j,B))
&\text{ in case }B\in R^{\flat}\\
\text{ord}_XH_0(X)+\sum_{1\le j\le\chi(H)}\text{ord}_XH_j(X,Y)
&\text{ in case }B\in R^{\flat}_{\infty}\\
\endcases
\tag{SP34}
$$
(with the understanding that if $H=0$ then $S(H,B)=\infty$), and we call 
$S(H,B)$ the {\bf $B$-strength} of $H$, and by (BP4) and (BP5)
we see that
$$
\cases
\text{for any $B\in R^{\flat}$ and $(z,V,W)\in\epsilon(B)$ we have}\\
\text{ord}_XH(X^V,z(X)+X^WY)=VS(H,B),\\
\text{and for any $B\in R^{\flat}_{\infty}$ we have}\\
\text{ord}_XH(X,Y)=S(H,B).\\
\endcases
\tag{SP35}
$$
We also put
$$
A^{**}(H,B)
=\prod_{1\le j\le\chi(H)\text{ with }H_j\in\tau^*(B)}A(R^*(H_j,B))
\quad\text{ for all }B\in R^{\flat *}
\tag{SP36}
$$
and we call $A^{**}(H,B)$ the {\bf doubly strict $B$-constant} of $H$, 
and we put
$$
D^{**}(H,B)=\left(\text{deg}_Y\Omega^*_B(H)\right)/D^*(B)
\quad\text{ for all }B\in R^{\flat *}
\tag{SP37}
$$
and we call $D^{**}(H,B)$ the {\bf doubly strict $B$-degree} of $H$, and we
note that, in view of (SBP4),
$$
\cases
\text{for any $B\in R^{\flat *}$ we have}\\
0\ne A^{**}(H,B)\in k\text{ and }0\le D^{**}(H,B)\in\Bbb Z\\
\text{with: }D^{**}(H,B)>0\Leftrightarrow\Omega^*_B(H)\ne 1
\endcases
\tag{SP38}
$$
and
$$
\cases
\text{for any $B\in R^{\flat *}$ and $(z,V,W)\in\epsilon(B)$ we have}\\
\text{minco}_X[z,V,W](\Omega^*_B(H))=A^{**}(H,B)E(B,Y)^{D^{**}(H,B)}.
\endcases
\tag{SP39}
$$

To collect together information about the $Y$-degrees of
$\Omega_B(H)$, $\Omega'_B(H)$, $\Omega^*_B(H)$, in view of (SP13) and
(SP27), we see that
$$
\cases
\text{for any $B\in R^{\flat}$ we have}\\
\text{deg}_Y\Omega_B(H)=\text{deg}_Y\Omega'_B(H)+\text{deg}_Y\Omega^*_B(H)\\
\text{and for any $(z,V,W)\in\epsilon(B)$ we have}\\
\text{deg}_Y\Omega'_B(H)
=D(B)\text{deg}_Y\text{pinco}_X[z,V,W](H,\sigma(B))\\
\text{and deg}_Y\Omega^*_B(H)
=D(B)\text{deg}_Y\text{sinco}_X[z,V,W](H,\sigma(B))\\
\text{and if }H\ne 0\text{ then we also have }\\
\text{deg}_Y\text{minco}_X[z,V,W](H)=
\text{deg}_Y\text{pinco}_X[z,V,W](H,\sigma(B))\\
\qquad\qquad\qquad\qquad\qquad\qquad\;
+~\text{deg}_Y\text{sinco}_X[z,V,W](H,\sigma(B))\\
\text{and deg}_Y\Omega_B(H)
=D(B)\text{deg}_Y\text{minco}_X[z,V,W](H).\\
\endcases
\tag{SP40}
$$
Next we recall that
$$
\cases
\text{for any $B\in\widehat R^{\flat}$ we have}\\
R^*(B)=\{B'\in R^{\flat *}\cup R^{\flat}_{\infty}:\lambda(B')
=\lambda(B)\text{ and }\sigma(B')=\tau^*(B')\cap\sigma(B)\}\\
\text{and }\sigma(B)=\coprod_{B'\in R^*(B)}\sigma(B')\\
\text{which is a partition of $\sigma(B)$
into pairwise disjoint nonempty subsets}\\
\endcases
\tag{SP41}
$$
and we put
$$
R^*(H,B)=\{B'\in R^*(B):\Omega^*_{B'}(H)\ne 1\}
\quad\text{ for all }B\in\widehat R^{\flat}
\tag{SP42}
$$
where we note that $R^*(H,B)$ is a finite set whose members may be called
the {\bf strict $B$-friends} of $H$. 
We also put
$$
D'(B)=\cases 
-D(B)+\sum_{B'\in R^*(B)}D^*(B')&\text{ for all }B\in R^{\flat}\\
-1+\sum_{B'\in R^*(B)}1&\text{ for all }B\in R^{\flat}_{\infty}\\
\endcases
\tag{SP43}
$$
(with the understanding that if $R^*(B)$ is an infinite set then 
$D'(B)=\infty$), and we call $D'(B)$ the {\bf primitive degree} of $B$.
\footnote{Observe that if $R^*(B)$ is a finite set then $D'(B)$ is a
nonnegative integer; in particular, if $B\in R^{\flat}_{\infty}$ then 
$\text{card}(R^*(B))=1$, where card denotes {\bf cardinality},
and hence $D'(B)=0$.
We shall use $D'(B)$ mainly when $R^*(H,B)=R^*(B)$; in 
that situation, $R^*(B)$ is obviously a finite set and hence $D'(B)$ is a 
nonnegative integer. For the use of $D'(B)$ 
in such a situation, see (SP60) and (SP70).}
Moreover we put
$$
D''(B)=\cases
-D(B)+\sum_{f\in\sigma(B)}\text{deg}_Yf&\text{ for all }B\in R^{\flat}\\
-1+\sum_{f\in\sigma(B)}\text{deg}_Yf&\text{ for all }B\in R^{\flat}_{\infty}\\
\endcases
\tag{SP44}
$$
(with the understanding that if $\sigma(B)$ is an infinite set then 
$D''(B)=\infty$), 
and we call $D''(B)$ the {\bf doubly primitive degree} of $B$,
\footnote{Observe that if $\sigma(B)$ is a finite set then $D''(B)$ is a
nonnegative integer. Also observe that Definitions (SP43) and (SP44)
would look more natural if for every $B\in R^{\flat}_{\infty}$ we put
$D(B)=D^*(B)=1$.}
and we note that
\footnote{An element $\Phi$ in $R$ is 
{\bf devoid of multiple factors} in $R$ means that $\Phi\ne 0$ and
the ideal $\Phi R$ is its own radical in $R$.}
$$
\cases
\text{if $B\in \widehat R^{\flat}$ is such that $R^*(H,B)=R^*(B)$}\\
\text{and $\Omega_B(H)$ is devoid of multiple factors in $R$,}\\
\text{then deg}_Y\Omega_B(H)=\cases
D(B)+D''(B)&\text{ in case }B\in R^{\flat}\\
1+D''(B)&\text{ in case }B\in R^{\flat}_{\infty}.\\
\endcases
\endcases
\tag{SP45}
$$
Now clearly
$$
\cases
\text{for all }B\in\widehat R^{\flat}\\
\text{we have }\Omega_B(H)=\Omega'_{B}(H)\Omega^*_{B}(H)\\
\text{and }\Omega^*_B(H)=\prod_{B'\in R^*(H,B)}\Omega^*_{B'}(H)
\endcases
\tag{SP46}
$$
and, in view of (SP23) to (SP26),
$$
\cases
\text{for all $B'\in R^*(B)$ with $B\in R^{\flat}$
we know that}\\
E(B',Y)\in k[Y]\text{ is monic in $Y$ having no multiple root in }k\\
\text{and deg}_YE(B',Y)=D^*(B')/D(B')=D^*(B')/D(B)>0
\endcases
\tag{SP47}
$$
and
$$
\cases
\text{for all $B'\ne B''$ in $R^*(B)$ with $B\in R^{\flat}$
we know that}\\
E(B',Y)\text{ and }E(B'',Y)\text{ have no common root in }k
\endcases
\tag{SP48}
$$
and
$$
\cases
\text{for all $B\in R^{\flat}$ with $(z,V,W)\in\epsilon(B)$
we know that}\\
\text{minco}_X[z,V,W](\Omega'_{B}(H))\text{ is a nonzero member of }k[Y]\\
\text{which has no common root with }E(B',Y)\text{ in }k\text{ for any }
B'\in R^*(B).
\endcases
\tag{SP49}
$$

By (SP10) and (SP20) we see that
$$
\cases
\text{for any $T\in R^{\sharp}$ we have:}\\
H_{\infty}=\prod_{R_{\infty}(T)\ne B\in T}\Omega'_B(H)\Leftrightarrow
\text{deg}_YH_{\infty}
=\sum_{R_{\infty}(T)\ne B\in T}\text{deg}_Y\Omega'_B(H).
\endcases
\tag{SP50}
$$
In the next Section we shall show that (SP50) is applicable
when $H=F_Y$ where $F\in R$ is devoid of multiple factors. In the Section 
after that we shall apply (SP30) to the case when $H$ is the
jacobian $J(F,G)$ of $F$ and $G$ in $R$. 
To prepare for all this, {\bf until further notice}, given any 
$B\in R^{\flat}$ with $(z,V,W)\in\epsilon(B)$, for every
$\Phi=\Phi(X,Y)\in R$ let us put
$$
\cases
\widetilde\Phi=\widetilde\Phi(X,Y)=\Phi(X^V,z(X)+X^WY)\\
\text{and}\\
I(\Phi)=\text{minco}_X[z,V,W](\Phi)=\text{inco}_X\widetilde\Phi.
\endcases
\tag{SP51}
$$
Then by (SP13) we see that
$$
\cases
\text{deg}_Y\Omega_B(H)=D(B)\text{deg}_YI(\Omega_B(H))\\
\text{and }I(H)=\mu I(\Omega_B(H))\\
\text{where $\mu\in k$ is such that: $\mu=0\Leftrightarrow H=0$.}
\endcases
\tag{SP52}
$$
and by the chain rule for partial derivatives we see that
$$ 
\cases
\text{if $I(H)\not\in k$}\\
\text{then }I(H_Y)=(I(H))_Y.
\endcases
\tag{SP53}
$$
Also clearly
$$
\cases
\text{for any $\Psi(Y)\in k[Y]\setminus k$
we have}\\
0\ne\Psi_Y(Y)\in k[Y]\text{ with }
\text{deg}_Y\Psi_Y(Y)=-1+\text{deg}_Y\Psi(Y)
\endcases
\tag{SP54}
$$
and
$$
\cases
\text{for any $\Psi(Y)\in k[Y]$ and $\mu\in k$ and $0<\nu\in\Bbb Z$
we have:}\\
\Psi(Y)=(Y-\mu)^{\nu}\Psi'(Y)
\text{ where }\Psi'(Y)\in k[Y]\text{ with }\Psi'(\mu)\ne 0\\
\Rightarrow\Psi_Y(Y)=(Y-\mu)^{\nu-1}\Psi''(Y)
\text{ where }\Psi''(Y)\in k[Y]\text{ with }\Psi''(\mu)\ne 0.
\endcases
\tag{SP55}
$$
By (SP46) we see that
$$
I(\Omega_B(H))
=I(\Omega'_{B}(H))\prod_{B'\in R^*(H,B)}I(\Omega^*_{B'}(H)).
\tag{SP56}
$$
Moreover, by (SP38), (SP39) and (SP47),
$$
\cases
\text{for all $B'\in R^*(H,B)$ we know that}\\
I(\Omega^*_{B'}(H))=A^{**}(H,B')E(B',Y)^{D^{**}(H,B')}\\
\text{where $0\ne A^{**}(H,B')\in k$, $0<D^{**}(H,B')\in\Bbb Z$,}\\
\text{$E(B',Y)\in k[Y]$ is monic in $Y$ having no multiple roots in $k$,}\\
\text{and deg$_YE(B',Y)=D^*(B')/D(B)>0$.}
\endcases
\tag{SP57}
$$
Likewise, by (SP48) and (SP49),
$$
\cases
\text{for all $B'\ne B''$ in $R^*(H,B)$ we know that}\\
\text{$E(B',Y)$ and $E(B'',Y)$ have no common root in $k$,}\\
\text{and we also know that
$I(\Omega'_{B}(H))$ is a nonzero member of $k[Y]$}\\
\text{which has no common root with $\prod_{B'\in R^*(H,B)}E(B',Y)$
in $k$.}
\endcases
\tag{SP58}
$$
By (SP51) to (SP58), we conclude that
$$
\cases
\text{if }\Omega_B(H)\ne 1\\
\text{then deg}_YI(\Omega_B(H_Y))=-1+\text{deg}_YI(\Omega_B(H))\\
\text{and }
I(\Omega_B(H_Y))=L(Y)\prod_{B'\in R^*(H,B)}E(B',Y)^{D^{**}(H,B')-1}\\
\text{where $0\ne L(Y)\in k[Y]$ has no common root with }
\prod_{B'\in R^*(H,B)}E(B',Y)\text{ in }k\\
\text{and deg}_YL(Y)=
-1+\text{deg}_YI(\Omega'_B(H))+\sum_{B'\in R^*(H,B)}\text{deg}_YE(B',Y)\\
\text{with deg}_YE(B',Y)=D^*(B')/D(B)\text{ for all }B'\in R^*(H,B).
\endcases
\tag{SP59}
$$

Now if $R^*(H,B)=R^*(B)$ then clearly $\Omega_B(H)\ne 1=\Omega'_B(H)$
and hence in particular deg$_YI(\Omega'_B(H))=0$ and therefore
in the situation of (SP59) we get
$$
\text{deg}_YL(Y)=-1+\sum_{B'\in R^*(B)}[D^*(B')/D(B)]=D'(B)/D(B)
$$ 
and in view of (SP26) we have 
$\text{pinco}_X[z,V,W](H_Y,\sigma(B))=\mu L(Y)$ with
$0\ne\mu\in k$, and hence
$$
\text{deg}_Y\text{pinco}_X[z,V,W](H_Y,\sigma(B))
=\text{deg}_YL(Y)=D'(B)/D(B).
$$ 
Thus
$$
\cases
\text{if $R^*(H,B)=R^*(B)$,}\\
\text{then }\Omega_B(H)\ne 1\\
\text{and deg}_Y\text{pinco}_X[z,V,W](H_Y,\sigma(B))=D'(B)/D(B).\\
\endcases
\tag{SP60}
$$
In view of (SP40) and (SP51), by (SP59) we see that
$$
\cases
\text{if }\Omega_B(H)\ne 1\\
\text{then deg}_Y\text{minco}_X[z,V,W](H_Y)
=-1+\text{deg}_Y\text{minco}_X[z,V,W](H)\\
\text{and deg}_Y\Omega_B(H_Y)=-D(B)+\text{deg}_Y\Omega_B(H).\\
\endcases
\tag{SP61}
$$
In view of (SP51), by (SP56) to (SP59) we also see that
$$
D^{**}(H_Y,B')=-1+D^{**}(H,B')
\quad\text{ for all }B'\in R^*(H,B).
\tag{SP62}
$$
By (SP62) it follows that
$$
\cases
\text{if $B\in R^{\flat *}$ and $D^{**}(H,B)>0$}\\
\text{then }D^{**}(H_Y,B)=-1+D^{**}(H,B).
\endcases
\tag{SP63}
$$
In view of (SP40), by (SP60) and (SP61) we see that
$$
\cases
\text{if $R^*(H,B)=R^*(B)$,}\\
\text{then deg}_Y\text{sinco}_X[z,V,W](H_Y,\sigma(B)) =-1
-[D'(B)/D(B)]+\text{deg}_Y\text{minco}_X[z,V,W](H).
\endcases
\tag{SP64}
$$
In view of (SP40), by (SP60) and (SP64) we see that
$$
\cases
\text{if $R^*(H,B)=R^*(B)$,}\\
\text{then deg}_Y\Omega_B(H_Y)=-D(B)+\text{deg}_Y\Omega_B(H)\\
\text{and deg}_Y\Omega'_B(H_Y)=D'(B)\\
\text{and deg}_Y\Omega^*_B(H_Y)=-D(B)-D'(B)+\text{deg}_Y\Omega_B(H).
\endcases
\tag{SP65}
$$

Turning to the jacobian, upon letting 
$$
\Delta=J(H,G)\quad\text{ and }\quad
\overline\Delta=J(\widetilde H,\widetilde G)
\quad\text{ with }G\in R
$$
by the chain rule for jacobians we get 
$$
\overline\Delta(X,Y)=VX^{V+W-1}\widetilde\Delta(X,Y)
$$
and now, assuming that $\Omega_B(H)\ne 1=\Omega_B(G)$ and
$G\ne 0\ne S(G,B)$, we have
$$
\widetilde H(X,Y)=I(H)X^{VS(H,B)}
+(\text{terms of degree}>VS(H,B)\text{ in }X)
$$
with $I(H)\in k[Y]\setminus k$, and
$$
\widetilde G(X,Y)=I(G)X^{VS(G,B)}
+(\text{terms of degree}>VS(G,B)\text{ in }X)
$$
with $0\ne I(G)\in k$ and $VS(G,B)\ne 0$, and hence we get
$$
\aligned
\overline\Delta(X,Y)=&-VS(G,B)I(G)(I(H))_YX^{VS(H,B)+VS(G,B)-1}\\
&+(\text{terms of degree}\ge[VS(H,B)+VS(G,B)]\text{ in }X)
\endaligned
$$
and therefore by (SP52) we have
$$
I(J(H,G))=\mu I(H_Y)\quad\text{ with }0\ne\mu\in k
$$
and hence by (SP51) we get
$$
\text{minco}_X[z,V,W](J(H,G))=\mu~\text{minco}_X[z,V,W](H_Y)
\quad\text{ with }0\ne\mu\in k.
$$
Thus
$$
\cases
\text{if $\Omega_B(H)\ne 1$ and $0\ne G\in R$ is such that
$\Omega_B(G)=1$ and $S(G,B)\ne 0$}\\
\text{then minco}_X[z,V,W](J(H,G))=\mu~\text{minco}_X[z,V,W](H_Y)
\text{ with }0\ne\mu\in k\\
\text{and hence pinco}_X[z,V,W](J(H,G),\sigma(B))
=\text{pinco}_X[z,V,W](H_Y,\sigma(B))\\
\text{and sinco}_X[z,V,W](J(H,G),\sigma(B))
=\text{sinco}_X[z,V,W](H_Y,\sigma(B)).\\
\endcases
\tag{SP66}
$$
In view of (SP40), by (SP66) we see that
$$
\cases
\text{if $\Omega_B(H)\ne 1$ and $0\ne G\in R$ is such that
$\Omega_B(G)=1$ and $S(G,B)\ne 0$}\\
\text{then deg}_Y\text{minco}_X[z,V,W](J(H,G))
=\text{deg}_Y\text{minco}_X[z,V,W](H_Y)\\
\text{and deg}_Y\text{pinco}_X[z,V,W](J(H,G),\sigma(B))
=\text{deg}_Y\text{pinco}_X[z,V,W](H_Y,\sigma(B))\\
\text{and deg}_Y\text{sinco}_X[z,V,W](J(H,G),\sigma(B))
=\text{deg}_Y\text{sinco}_X[z,V,W](H_Y,\sigma(B))\\
\text{and deg}_Y\Omega_B(J(H,G))=\text{deg}_Y\Omega_B(H_Y)\\
\text{and deg}_Y\Omega'_B(J(H,G))=\text{deg}_Y\Omega'_B(H_Y)\\
\text{and deg}_Y\Omega^*_B(J(H,G))=\text{deg}_Y\Omega^*_B(H_Y).
\endcases
\tag{SP67}
$$

To consider another similar case, just for a moment let $j(H,G)$ stand for
$(HG)_Y$; then upon letting
$$
\delta=j(H,G)\quad\text{ and }\quad
\overline\delta=j(\widetilde H,\widetilde G)
\quad\text{ with }G\in R
$$
by the product rule for derivatives we get 
$$
\overline\delta(X,Y)=X^W\widetilde\delta(X,Y)
$$
and now, assuming that $\Omega_B(H)\ne 1=\Omega_B(G)$ and
$G\ne 0$, we have
$$
\widetilde H(X,Y)=I(H)X^{VS(H,B)}
+(\text{terms of degree}>VS(H,B)\text{ in }X)
$$
with $I(H)\in k[Y]\setminus k$, and
$$
\widetilde G(X,Y)=I(G)X^{VS(G,B)}
+(\text{terms of degree}>VS(G,B)\text{ in }X)
$$
with $0\ne I(G)\in k$, and hence we get
$$
\aligned
\overline\delta(X,Y)=&I(G)(I(H))_YX^{VS(H,B)+VS(G,B)}\\
&+(\text{ terms of degree}>[VS(H,B)+VS(G,B)]\text{ in }X)
\endaligned
$$
and therefore by (SP52) we have
$$
I((HG)_Y)=\mu I(H_Y)\quad\text{ with }0\ne\mu\in k
$$
and hence by (SP51) we get
$$
\text{minco}_X[z,V,W]((HG)_Y)=\mu~\text{minco}_X[z,V,W](H_Y)
\quad\text{ with }0\ne\mu\in k.
$$
Thus
$$
\cases
\text{if $\Omega_B(H)\ne 1$ and $0\ne G\in R$ is such that
$\Omega_B(G)=1$}\\
\text{then minco}_X[z,V,W]((HG)_Y)=\mu~\text{minco}_X[z,V,W](H_Y)
\text{ with }0\ne\mu\in k\\
\text{and hence pinco}_X[z,V,W]((HG)_Y,\sigma(B))
=\text{pinco}_X[z,V,W](H_Y,\sigma(B))\\
\text{and sinco}_X[z,V,W]((HG)_Y,\sigma(B))
=\text{sinco}_X[z,V,W](H_Y,\sigma(B)).\\
\endcases
\tag{SP68}
$$
In view of (SP40), by (SP68) we see that
$$
\cases
\text{if $\Omega_B(H)\ne 1$ and $0\ne G\in R$ is such that
$\Omega_B(G)=1$}\\
\text{then deg}_Y\text{minco}_X[z,V,W]((HG)_Y)
=\text{deg}_Y\text{minco}_X[z,V,W](H_Y)\\
\text{and deg}_Y\text{pinco}_X[z,V,W]((HG)_Y,\sigma(B))
=\text{deg}_Y\text{pinco}_X[z,V,W](H_Y,\sigma(B))\\
\text{and deg}_Y\text{sinco}_X[z,V,W]((HG)_Y,\sigma(B))
=\text{deg}_Y\text{sinco}_X[z,V,W](H_Y,\sigma(B))\\
\text{and deg}_Y\Omega_B((HG)_Y)=\text{deg}_Y\Omega_B(H_Y)\\
\text{and deg}_Y\Omega'_B((HG)_Y)=\text{deg}_Y\Omega'_B(H_Y)\\
\text{and deg}_Y\Omega^*_B((HG)_Y)=\text{deg}_Y\Omega^*_B(H_Y).
\endcases
\tag{SP69}
$$

{\bf Abandoning notation (SP51),} let us summarize Results 
(SP60) to (SP69) as Lemmas (SP70) to (SP75) stated below.

\proclaim{Lemma (SP70)}
If $B\in R^{\flat}$ and $F\in R$ are such that $R^*(F,B)=R^*(B)$,
then $\Omega_B(F)\ne 1$ and we have
$$
\cases
\text{deg}_Y\Omega_B(F_Y)=-D(B)+\text{deg}_Y\Omega_B(F)\\
\text{and deg}_Y\Omega'_B(F_Y)=D'(B)\\
\text{and deg}_Y\Omega^*_B(F_Y)=-D(B)-D'(B)+\text{deg}_Y\Omega_B(F)
\endcases
$$
and for every $(z,V,W)\in\epsilon(B)$ we have
$$
\cases
\text{deg}_Y\text{minco}_X[z,V,W](F_Y)
=-1+\text{deg}_Y\text{minco}_X[z,V,W](F)\\
\text{and deg}_Y\text{pinco}_X[z,V,W](F_Y,\sigma(B))=D'(B)/D(B)\\
\text{and deg}_Y\text{sinco}_X[z,V,W](F_Y,\sigma(B)) =-1
-[D'(B)/D(B)]+\text{deg}_Y\text{minco}_X[z,V,W](F).
\endcases
$$
\endproclaim

\proclaim{Lemma (SP71)}
If $B\in R^{\flat}$ and $F\in R$ are such that $\Omega_B(F)\ne 1$,
then we have
$$
\text{deg}_Y\Omega_B(F_Y)=
-D(B)+\text{deg}_Y\Omega_B(F)
$$
and for every $(z,V,W)\in\epsilon(B)$ we have
$$
\text{deg}_Y\text{minco}_X[z,V,W](F_Y)
=-1+\text{deg}_Y\text{minco}_X[z,V,W](F).
$$
\endproclaim

\proclaim{Lemma (SP72)}
Given any $B\in R^{\flat}$ and $F\in R$, for every $B'\in R^*(F,B)$ we have
$$
D^{**}(F_Y,B)=-1+D^{**}(F,B).
$$
\endproclaim

\proclaim{Lemma (SP73)}
If $B\in R^{\flat *}$ and $F\in R$ are such that $D^{**}(F,B)>0$,
then we have
$$
D^{**}(F_Y,B)=-1+D^{**}(F,B).
$$
\endproclaim

\proclaim{Lemma (SP74)}
If $B\in R^{\flat}$ and $F\in R$ are such that $\Omega_B(F)\ne 1$, and 
$0\ne G\in R$ is such that $\Omega_B(G)=1$ and $S(G,B)\ne 0$,
then we have
$$
\cases
\text{deg}_Y\Omega_B(J(F,G))=\text{deg}_Y\Omega_B(F_Y)\\
\text{and deg}_Y\Omega'_B(J(F,G))=\text{deg}_Y\Omega'_B(F_Y)\\
\text{and deg}_Y\Omega^*_B(J(F,G))=\text{deg}_Y\Omega^*_B(F_Y)
\endcases
$$
and for every $(z,V,W)\in\epsilon(B)$ we have
$$
\cases
\text{deg}_Y\text{minco}_X[z,V,W](J(F,G),\sigma(B))
=\text{deg}_Y\text{minco}_X[z,V,W](F_Y,\sigma(B))\\
\text{and deg}_Y\text{pinco}_X[z,V,W](J(F,G),\sigma(B))
=\text{deg}_Y\text{pinco}_X[z,V,W](F_Y,\sigma(B))\\
\text{and deg}_Y\text{sinco}_X[z,V,W](J(F,G),\sigma(B))
=\text{deg}_Y\text{sinco}_X[z,V,W](F_Y,\sigma(B))
\endcases
$$
and actually we have
$$
\cases
\text{minco}_X[z,V,W](J(F,G))=\mu~\text{minco}_X[z,V,W](F_Y)
\text{ with }0\ne\mu\in k,\\
\text{and pinco}_X[z,V,W](J(F,G),\sigma(B))
=\text{pinco}_X[z,V,W](F_Y,\sigma(B)),\\
\text{and sinco}_X[z,V,W](J(F,G),\sigma(B))
=\text{sinco}_X[z,V,W](F_Y,\sigma(B)).\\
\endcases
$$
\endproclaim

\proclaim{Lemma (SP75)}
If $B\in R^{\flat}$ and $F\in R$ are such that $\Omega_B(F)\ne 1$, and 
$0\ne G\in R$ is such that $\Omega_B(G)=1$,
then we have
$$
\cases
\text{deg}_Y\Omega_B((FG)_Y)=\text{deg}_Y\Omega_B(F_Y)\\
\text{and deg}_Y\Omega'_B((FG)_Y)=\text{deg}_Y\Omega'_B(F_Y)\\
\text{and deg}_Y\Omega^*_B((FG)_Y)=\text{deg}_Y\Omega^*_B(F_Y)
\endcases
$$
and for every $(z,V,W)\in\epsilon(B)$ we have
$$
\cases
\text{deg}_Y\text{minco}_X[z,V,W]((FG)_Y)
=\text{deg}_Y\text{minco}_X[z,V,W](F_Y)\\
\text{and deg}_Y\text{pinco}_X[z,V,W]((FG)_Y,\sigma(B))
=\text{deg}_Y\text{pinco}_X[z,V,W](F_Y,\sigma(B))\\
\text{and deg}_Y\text{sinco}_X[z,V,W]((FG)_Y,\sigma(B))
=\text{deg}_Y\text{sinco}_X[z,V,W](F_Y,\sigma(B))\\
\endcases
$$
and actually we have
$$
\cases
\text{minco}_X[z,V,W]((FG)_Y)=\mu~\text{minco}_X[z,V,W](F_Y)
\text{ with }0\ne\mu\in k,\\
\text{and pinco}_X[z,V,W]((FG)_Y,\sigma(B))
=\text{pinco}_X[z,V,W](F_Y,\sigma(B)),\\
\text{and sinco}_X[z,V,W]((FG)_Y,\sigma(B))
=\text{sinco}_X[z,V,W](F_Y,\sigma(B)).\\
\endcases
$$
\endproclaim

Now, as a consequence of (SP45) and (SP70) we shall 
prove the following Lemma:

\proclaim{Lemma (SP76)}
If $B\in\widehat R^{\flat}$ and $F\in R$ are such that $R^*(F,B)=R^*(B)$
and $F$ is devoid of multiple factors in $R$, then
$\Omega_B(F)\ne 1$ and
$\text{deg}_Y\Omega_B(F_Y)=D''(B)$ and $\text{deg}_Y\Omega'_B(F_Y)=D'(B)$.
\endproclaim

Namely, if $B\in R^{\flat}$ and $F\in R$ are such that $R^*(F,B)=R^*(B)$
then by (SP70) we get $\Omega_B(F)\ne 1$ and
$\text{deg}_Y\Omega_B(F_Y)=-D(B)+\text{deg}_Y\Omega_B(F)$ and 
$\text{deg}_Y\Omega'_B(F_Y)=D'(B)$; if $F$ is also devoid of multiple factors 
in $R$, then by (SP45) we know that
$\text{deg}_Y\Omega_B(F)=D(B)+D''(B)$ and and hence we get
$\text{deg}_Y\Omega_B(F_Y)=D''(B)$.
Likewise, if $B\in R^{\flat}_{\infty}$ and $F\in R$ are such that
$\Omega_B(F)\ne 1$ then clearly 
$\text{deg}_Y\Omega_B(F_Y)=-1+\text{deg}_Y\Omega_B(F)$; 
if $F$ is also devoid of multiple factors in $R$, then by (SP45) we know that
$\text{deg}_Y\Omega_B(F)=1+D''(B)$ and and hence we get
$\text{deg}_Y\Omega_B(F_Y)=D''(B)$. 
This completes the proof of (SP76).
\footnote{We may tacitly use the {\bf obvious facts} that:
(1) if $B\in\widehat R^{\flat}$ and $F\in R$ are such that 
$R^*(F,B)=R^*(B)$ then $\Omega_B(F)\ne 1$; (2) for every 
$B\in R^{\flat}_{\infty}$ we have $D'(B)=0$; (3) for every
$B\in R^{\flat}_{\infty}$ and every $G\in R$ we have $\Omega_B(G)=1$
and hence $\text{deg}_Y\Omega'_B(G)=0$.}

\centerline{}

Next, as a consequence of (SP74) and (SP75) we shall prove the following 
Lemma:

\proclaim{Lemma (SP77)}
Given any $F\in R\setminus k((X)$ and $0\ne G\in R$, upon letting $T=T(FG)$, 
we have the following.

(SP77.1) If $B\in T$ is such that $\Omega_B(G)=1$ 
and $S(G,B)\ne 0$ then
$$
\Omega_B(F)\ne 1
$$
and
$$
\text{deg}_Y\Omega_B(J(F,G))=\text{deg}_Y\Omega_B(F_Y)
\quad\text{ and }\quad
\text{deg}_Y\Omega'_B(J(F,G))=\text{deg}_Y\Omega'_B(F_Y)
$$
and for every $B'\in\pi(T,B)$ we have 
$$
\Omega_{B'}(F)\ne 1=\Omega_{B'}(G)
\quad\text{ and }\quad
S(G,B')\ne 0
$$
and
$$
\text{deg}_Y\Omega_{B'}(J(F,G))=\text{deg}_Y\Omega_{B'}(F_Y)
\quad\text{ and }\quad
\text{deg}_Y\Omega'_{B'}(J(F,G))=\text{deg}_Y\Omega'_{B'}(F_Y).
$$

(SP77.2) If $B\in T$ is such that $\Omega_B(G)=1$ then
$$
\Omega_B(F)\ne 1
$$
and
$$
\text{deg}_Y\Omega_B((FG)_Y)=\text{deg}_Y\Omega_B(F_Y)
\quad\text{ and }\quad
\text{deg}_Y\Omega'_B((FG)_Y)=\text{deg}_Y\Omega'_B(F_Y)
$$
and for every $B'\in\pi(T,B)$ we have 
$$
\Omega_{B'}(F)\ne 1=\Omega_{B'}(G)
$$
and 
$$
\text{deg}_Y\Omega_{B'}((FG)_Y)=\text{deg}_Y\Omega_{B'}(F_Y)
\quad\text{ and }\quad
\text{deg}_Y\Omega'_{B'}((FG)_Y)=\text{deg}_Y\Omega'_{B'}(F_Y).
$$
\endproclaim

Namely, if $B\in T$ is such that $\Omega_B(G)=1$ then obviously
$\Omega_B(F)\ne 1$ and for
every $B'\in\pi(T,B)$ we have $\Omega_{B'}(F)\ne 1=\Omega_{B'}$, and
by (BP5) we also see that if $B\in T$ is such that $\Omega_B(G)=1$ then 
for every $B'\in\pi(T,B)$ we have $S(G,B')=S(G,B)$. Therefore, in case of 
$B\in T\setminus\{R_{\infty}(T)\}$, our assertions follow from
(SP74) and (SP75). 
Moreover, if $B=R_{\infty}(T)$ and $\Omega_B(G)=1$ and
$S(G,B)\ne 0$, then by the equation $J(F,G)=F_XG_Y-F_YG_X$
we see that $J(F,G)=-F_YG_X$ with $0\ne G_X\in k((X))$ and hence
$\text{deg}_Y\Omega_B(J(F,G))=\text{deg}_Y\Omega_B(F_Y)$ and
$\text{deg}_Y\Omega'_B(J(F,G))=0=\text{deg}_Y\Omega'_B(F_Y)$.
Likewise, if $B=R_{\infty}(T)$ and $\Omega_B(G)=1$,
then by the equation $(FG)_Y=FG_Y+F_YG$
we see that $(FG)_Y=F_YG$ with $0\ne G\in k((X))$ and hence
$\text{deg}_Y\Omega_B((FG)_Y)=\text{deg}_Y\Omega_B(F_Y)$ and
$\text{deg}_Y\Omega'_B((FG)_Y)=0=\text{deg}_Y\Omega'_B(F_Y)$.
This completes the proof of (SP77).

\centerline{}

Now we shall prove the following Lemma:

\proclaim{Lemma (SP78)}
For any $\Phi\in R$, upon letting $T=T(\Phi)$, we have the following.

(SP78.1) Given any $B\in T$ with $\pi(T,B)\ne\emptyset$, 
for every $B'\in\rho(T,B)$ there is a unique $\alpha(B')\in R^*(B)$ with
$\sigma(\alpha(B'))=\sigma(B')$, and $B'\mapsto\alpha(B')$ gives a bijection
of $\rho(T,B)$ onto $R^*(B)$. Moreover, for every $B'\in\rho(T,B)$ we have 
$$
D(B')=\cases D^*(\alpha(B'))&\text{ in case }B\in R^{\flat}\\
1&\text{ in case }B\in R^{\flat}_{\infty}\\
\endcases
$$ 
and hence we have 
$$
D'(B)=\cases
-D(B)+\sum_{B'\in\rho(T,B)}D(B')&\text{ in case }B\in R^{\flat}\\
-1+\sum_{B'\in\rho(T,B)}D(B')&\text{ in case }B\in R^{\flat}_{\infty}
\endcases
$$
where we note that if $B\in R^{\flat}_{\infty}$ then
$\text{card}(\rho(T,B))=\text{card}(R^*(B))=1$. 

(SP78.2) Given any $B\in T$ with $\pi(T,B)=\emptyset$, 
for every $f\in\sigma(B)$ there is a unique $\beta(f)\in R^*(B)$
with $f\in\sigma(\beta(f))$, and $f\mapsto\beta(f)$ gives a bijection
of $\sigma(B)$ onto $R^*(B)$. Moreover, for every $f\in\sigma(B)$ we have 
$$
\text{deg}_Yf=\cases
D^*(\beta(f))&\text{ in case }B\in R^{\flat}\\
1&\text{ in case }B\in R^{\flat}_{\infty}
\endcases
$$ 
and hence we have 
$$
D''(B)=D'(B).
$$

(SP78.3) For every $B\in T$ we have
$$
D''(B)=D'(B)+\sum_{B'\in\pi(T,B)}D'(B').
$$
\endproclaim

Namely, the proofs of (SP78.1) and (SP78.2) are straightforward.
We shall prove (SP78.3) by induction on card$(\pi(T,B))$.
In case of card$(\pi(T,B))=0$ our assertion follows from (SP78.2).
So let card$(\pi(T,B))>0$ and assume true for all
smaller values of card$(\pi(T,B))$. Then, by the induction hypothesis,
for every $B'\in\rho(T,B)$ we have
$$
D''(B')=D'(B')+\sum_{B''\in\pi(T,B')}D'(B'')
$$
and hence by the definition of $D''(B')$ we get
$$
D'(B')+\sum_{B''\in\pi(T,B')}D'(B'')
=-D(B')+\sum_{f\in\sigma(B')}\text{deg}_Yf.
$$
Summing both sides of the above equation as $B'$ varies over $\rho(T,B)$,
we get
$$
\sum_{B'\in\pi(T,B)}D'(B')
=-\left[\sum_{B'\in\rho(T,B)}D(B')\right]
+\left[\sum_{f\in\sigma(B)}\text{deg}_Yf\right]
$$
and by (SP78.1) we have
$$
D'(B)=\cases
-D(B)+\sum_{B'\in\rho(T,B)}D(B')&\text{ in case }B\in R^{\flat}\\
-1+\sum_{B'\in\rho(T,B)}D(B')&\text{ in case }B\in R^{\flat}_{\infty}.\\
\endcases
$$
By adding the above two equations we get
$$
D'(B)+\sum_{B'\in\pi(T,B)}D'(B')=\cases
-D(B)+\sum_{f\in\sigma(B)}\text{deg}_Yf
&\text{ in case }B\in R^{\flat}\\
-1+\sum_{f\in\sigma(B)}\text{deg}_Yf
&\text{ in case }B\in R^{\flat}_{\infty}\\
\endcases
$$
and hence by the definition of $D''(B)$ we conclude that
$$
D''(B)=D'(B)+\sum_{B'\in\pi(T,B)}D'(B').
$$
This completes the proof of (SP78.3).

\centerline{}

As an immediate consequence of (SP78.3) we get the following Lemma:

\proclaim{Lemma (SP79)}
Let $B\in T=T(\Phi)$ with $\Phi\in R$ be such that: 
$\text{deg}_Y\Omega_B(H)=D''(B)$, $\text{deg}_Y\Omega'_B(H)=D'(B)$, and 
$\text{deg}_Y\Omega'_{B'}(H)=D'(B')$
for all $B'\in\pi(T,B)$. Then
$$
\text{deg}_Y\Omega_B(H)
=\text{deg}_Y\Omega'_B(H)+\sum_{B'\in\pi(T,B)}\text{deg}_Y\Omega'_{B'}(H).
$$
\endproclaim

For any $B\in T=T(\Phi)$ with $\Phi\in R$, it is clear that
$\Omega'_B(H)\prod_{B'\in\pi(T,B)}\Omega'_{B'}(H)$ divides $\Omega_B(H)$
in $R$ and hence, as a companion to (SP50), and as a principle
applicable in the situation of (SP79), we get the following Lemma:

\proclaim{Lemma (SP80)}
For any $B\in T=T(\Phi)$ with $\Phi\in R$ we have:
$$
\cases
\Omega_B(H)=\Omega'_B(H)\prod_{B'\in\pi(T,B)}\Omega'_{B'}(H)\\
\Leftrightarrow\text{deg}_Y\Omega_B(H)=\text{deg}_Y\Omega'_B(H)
+\sum_{B'\in\pi(T,B)}\text{deg}_Y\Omega'_{B'}(H).
\endcases
$$
\endproclaim

\centerline{}

\centerline{{\bf Section 8: Factorization of the Derivative}}

If $T=T(F)$ where $F\in R\setminus k((X))$, then for every $B\in T$
we clearly have $R^*(F,B)=R^*(B)$. Therefore by (SP76), (SP79) and (SP80) 
we get the following Derivative Factorization Theorem.

\proclaim{Theorem (DF1)}
Let $T=T(F)$ where $F\in R\setminus k((X))$ is devoid of multiple 
factors in $R$. Then we have the following. 

(DF1.1) For any $B\in T$ we have
$$
\text{deg}_Y\Omega_B(F_Y)=D''(B)
\quad\text{ and }\quad
\text{deg}_Y\Omega'_B(F_Y)=D'(B)
$$
and
$$
\Omega_B(F_Y)=\Omega'_B(F_Y)\prod_{B'\in\pi(T,B)}\Omega'_{B'}(F_Y).
$$
where for every $B'\in\pi(T,B)$ we have
$$
\text{deg}_Y\Omega_{B'}(F_Y)=D''(B')
\quad\text{ and }\quad
\text{deg}_Y\Omega'_{B'}(F_Y)=D'(B').
$$

(DF1.2) By taking $B=R_{\infty}(T)$ in (DF1.1), 
for the monic part $(F_Y)_{\infty}$ of $F_Y$ we get
$$
(F_Y)_{\infty}=\prod_{B\in T\setminus\{R_{\infty}(T)\}}\Omega'_{B}(F_Y).
$$
\endproclaim

\centerline{}

\noindent
{\bf Remark (DF2).}
In the Factorization (DF1.2), 
the factor $\Omega'_B(F_Y)$ {\bf really occurs}, i.e., its $Y$-degree
$D'(B)$ is nonzero, if and only if either: (*) card$(R^*(B))=1$ and for the
unique $B'\in R^*(B)$ we have $D^*(B')>D(B)$, or: (**) card$(R^*(B))>1$.
Note that in the irreducible case, i.e., when $F=f\in R^{\natural}$,
(*) is always satisfied. Moreover, in the {\bf nontrivial irreducible case}, 
i.e., when $F=f\in R^{\natural}$ with deg$_Yf=n>1$, let us put
$$
\widehat h=\cases
h(c(f))&\text{ if }c_1(f)\not\in\Bbb Z\\
h(c(f))-1&\text{ if }c_1(f)\in\Bbb Z\\
\endcases
$$
and for $1\le i\le\widehat h$ let us put
$$
\widehat c_i=\cases
c_i(f)&\text{ if }c_1(f)\not\in\Bbb Z\\
c_{i+1}(f)&\text{ if }c_1(f)\in\Bbb Z\\
\endcases
$$
and
$$
\widehat r_i=\cases
r_i(q(m(f)))&\text{ if }r_1(f)\not\in\Bbb Z\\
r_{i+1}(q(m(f)))&\text{ if }r_1(f)\in\Bbb Z\\
\endcases
$$
and for $1\le i\le\widehat h+1$ let us put
$$
\widehat d_i=\cases
d_i(m(f))&\text{ if }c_1(f)\not\in\Bbb Z\\
d_{i+1}(m(f))&\text{ if }c_1(f)\in\Bbb Z.\\
\endcases
$$
Then $\widehat h$ is a positive 
integer, $\widehat c_1<\widehat c_2<\dots<\widehat c_{\widehat h}$ are in
$\Bbb Q\setminus\Bbb Z$, and 
$n=\widehat d_1>\widehat d_2>\dots>\widehat d_{\widehat h+1}=1$ are 
integers with $\widehat d_i\equiv 0\pmod{\widehat d_{i+1}}$ for 
$1\le i\le \widehat h$.  
Let $B_0=(\sigma(B_0),\lambda(B_0))\in R_{\infty}^{\flat}$
with $\sigma(B_0)=\{f\}$ and $\lambda(B_0)=-\infty$.
For $1\le i\le \widehat h$ let $B_i=(\sigma(B_i),\lambda(B_i))\in R^{\flat}$
with $\sigma(B_i)=\{f\}$ and $\lambda(B_i)=\widehat c_i$. Then
$T=T(f)=\{B_0,B_1,\dots,B_{\widehat h}\}$ with 
$R_{\infty}(T)=B_0<B_1<\dots<B_{\widehat h}$, and
for $1\le i\le \widehat h$ we have 
$$
S(B_i)=(\widehat d_i\widehat r_i)/n^2
\quad\text{ and }\quad
D(B_i)=n/\widehat d_i
$$ 
and
$$
D^*(B'_i)=n/\widehat d_{i+1}\quad\text{ where }\quad\{B'_i\}=R^*(B_i).
$$ 
Now
$$
(f_Y)_{\infty}=(1/n)f_Y=\prod_{1\le i\le\widehat h}\Omega'_{B_i}(f_Y)
$$
and for $1\le i\le\widehat h$ we have
$$
\text{deg}_Y\Omega'_{B_i}(f_Y)=D'(B_i)=-D(B_i)+D^*(B'_i) 
=(n/\widehat d_{i+1})-(n/\widehat d_i)>0.
$$
Let us factor $f_Y$ into irreducible factors by writing
$$
f_Y=n\prod_{1\le j\le\chi}f^{(j)}\quad\text{ with }f^{(j)}\in R^{\natural}
$$
and for $1\le i\le\widehat h$ let us put
$$
i^*=\{j\in\{1,\dots,\chi\}:\text{noc}(f,f^{(j)})=\widehat c_i\}.
$$
Then
$$
\{1,\dots,\chi\}=\coprod_{1\le i\le\widehat h}i^*
$$
is a partition into pairwise disjoint nonempty sets, and for
$1\le i\le\widehat h$ we have
$$
\Omega'_{B_i}(f_Y)=\prod_{j\in i^*}f^{(j)}
\quad\text{ with }\quad 0<\text{deg}_Yf^{(j)}\in(n/\widehat d_i)\Bbb Z
\quad\text{ for all }j\in i^*
$$
and
$$
\text{int}(f,\Omega'_{B_i}(f_Y))
=nS(B_i)\text{deg}_Y\Omega'_{B_i}(f_Y)
=[(\widehat d_i/\widehat d_{i+1})-1]\widehat r_i
$$ 
where int denotes intersection multiplicity.
\footnote{The {\bf intersection multiplicity} int$(f,g)$
of $f\in R^{\natural}$ with $g\in R$ is defined by putting
int$(f,g)=\text{ord}_X\text{Res}_Y(f,g)$, where Res$_Y(f,g)$ denotes the
$Y$-resultant of $f$ and $g$; equivalently, for any $z(X)\in k((X))$
with $f(X^n,z(X))=0$ where deg$_Yf=n$, we have
int$(f,g)=\text{ord}_Xg(X^n,z(X))$; see pages 286-287 of \cite{Ab}.
By (GNP7) we see that, given any $B\in R^{\flat}$ and $H\in R$,
for every $f\in\sigma(B)$ we have
int$(f,\Omega'_B(H))=nS(B)\text{deg}_Y\Omega'_B(H)$ where
deg$_Yf=n$.}

\centerline{}

\noindent
{\bf Example (DF3).} Now, if we are in the nontrivial irreducible case of
$F=f\in R^{\natural}$ with deg$_Yf=n>1$, and if $h(c(f))=1$ with
$c_1(f)\not\in\Bbb Z$, then the conclusions of the above Remark (DF2) say
that $\Omega'_{B_1}(f_Y)=(1/n)f_Y$ with int$(f,f_Y)=(n-1)m_1(f)$, and
noc$(f,f^{(j)})=c_1(f)=m_1(f)/n$ 
for every irreducible factor $f^{(j)}$ of $f_Y$.
To verify this in a particular situation,
by taking $(w_1(X),\dots,w_{n-1}(X),w_n(X))=(0,\dots,0,X^e)$
in Example (TR3) of Section 6 we have
$$
F(X,Y)=f(X,Y)=Y^n+X^e\in R^{\natural}
\text{ where }0\ne e\in\Bbb Z\text{ with GCD}(n,e)=1
$$
and hence $h(T(f))=h(f)=1$ with $m_1(f)=e$ and 
$$
l_0(T(f))=-\infty\text{ and }l_1(T(f))=c_1(f)=e/n
$$
and upon letting 
$$
B_i\in\widehat R^{\flat}\text{ with }
\sigma(B_i)=\{f\}\text{ and }\lambda(B_i)=l_i(T(f))
\text{ for }0\le i\le 1
$$
we have
$$
T(F)=T(f)=\{B_0,B_1\}
$$
with
$$
D'(B_0)=0\text{ and }D'(B_1)=n-1.
$$
Now clearly $f_Y=nY^{n-1}$, and hence
Res$_Y(f,f_Y)=n^n X^{(n-1)e}$ and $f^{(j)}=Y$ for
$1\le j\le\chi=n-1$, and therefore
int$(f,f_Y)=\text{ord}_X\text{Res}_Y(f,f_Y)=(n-1)e=(n-1)m_1(f)$ and
noc$(f,f^{(j)})=(1/n)\text{ord}_Xf(X,0)=m_1(f)/n$
for $1\le j\le\chi=n-1$. This completes the verification.

\centerline{}

\noindent
{\bf Example (DF4).} Next, if we are in the nontrivial irreducible case of
$F=f\in R^{\natural}$ with deg$_Yf=n>1$, and if $h(c(f))=2$ with
$c_1(f)\not\in\Bbb Z$, then the conclusions of the above Remark (DF2) say
that $\Omega'_{B_1}(f_Y)\Omega'_{B_2}(f_Y)=(1/n)f_Y$ and for
$1\le i\le 2$ we have that: 
deg$_Y\Omega'_{B_i}(f_Y)=D'(B_i)=(n/d_{i+1})-(n/d_i)>0$ and
int$(f,\Omega'_{B_i}(f_Y))=[(d_i/d_{i+1})-1]r_i$ where $\Omega'_{B_i}(f_Y)$ 
is the product of all those irreducible factors $f^{(j)}$ of $f_Y$ for which
noc$(f,f^{(j)})=c_i(f)$, and moreover for each of these $f^{(j)}$ we have
$0<\text{deg}_Yf^{(j)}\in(n/d_i)\Bbb Z$.
To verify this in a particular situation, 
in Example (TR4) of Section 6 we have
$$
F(X,Y)=f(X,Y)=(Y^2-X^{2a+1})^2-X^{3a+b+2}Y\in R^{\natural}
\text{ with }a\in\Bbb Z\text{ and }0\le b\in\Bbb Z
$$
with
$$
\cases
n=d_1=4\text{ and }d_2=2\text{ and }d_3=1,\text{ and }\\
c_1(f)=(2a+1)/2\text{ and }c_2(f)=(4a+2b+3)/4,\text{ and }\\
[(d_1/d_2)-1]r_1=4a+2\text{ and }[(d_2/d_3)-1]r_2=8a+2b+5
\endcases
$$
and $h(T(f))=2$ with $l_0(T(f))=-\infty$ and
$$
l_1(T(f))=c_1(f)=(2a+1)/2\text{ and }l_2(T(f))=c_2(f)=(4a+2b+3)/4
$$
and upon letting 
$$
B_i\in\widehat R^{\flat}\text{ with }
\sigma(B_i)=\{f\}\text{ and }\lambda(B_i)=l_i(T(f))
\text{ for }0\le i\le 2
$$
we have
$$
T(F)=T(f)=\{B_0,B_1,B_2\}.
$$
with $D'(B_0)=0$ and 
$$
D'(B_1)=1\text{ and }D'(B_2)=2.
$$

Now 
$$
f_Y=4Y(Y^2-X^{2a+1})-X^{3a+b+2}
$$
and hence by (TR5) of Section 6 we see that
$f_Y=4Yf^{(1)}f^{(2)}$ where
$$
f^{(1)}(X,Y)=Y-v(X)\in R^{\natural}\text{ and }v(X)\in k((X))
$$
with
$$
\text{ord}_Xv(X)=a+b+1
$$
and
$$
f^{(2)}(X,Y)=Y^2+\sum_{1\le i\le 2}v'_i(X)Y^{2-i}\in R^{\natural}
\text{ and }v'_i(X)\in k((X))
$$
with
$$
\text{ord}_Xv'_1(X)>(2a+1)/2\text{ and }\text{ord}_Xv'_2(X)=2a+1.
$$
Comparing coefficients of $Y^2$ and $Y$ in the equation 
$f_Y=4f^{(1)}f^{(2)}$ we see that 
$v'_1(X)-v(X)=0$ and $v'_2(X)-v'_1(X)v(X)=-X^{2a+1}$, and hence
$$
f^{(2)}(X,Y)=Y^2+v(X)Y-X^{2a+1}+v(X)^2.
$$
Applying the quadratic equation formula to the above equation we get the
roots of $f^{(2)}(X^4,Y)$ to be
$$
\align
Y&=(-1/2)v(X^4)\pm(1/2)\sqrt{4X^{8a+4}-3v(X^4)^2}\\
&=(-1/2)v(X^4)\pm X^{4a+2}\sqrt{1-(3/4)X^{-8a-4}v(X^4)^2}\\
&=(-1/2)[\mu X^{4a+4b+4}+(\text{terms of degree $>4a+4b+4$ in $X$})]\\
&\quad\pm X^{4a+2}
\sqrt{1-(3/4)\mu^2 X^{8b+4}+(\text{terms of degree $>8b+4$ in $X$})}\\
&=(-1/2)[\mu X^{4a+4b+4}+(\text{terms of degree $>4a+4b+4$ in $X$})]\\
&\quad\pm X^{4a+2}
[1-(3/8)\mu^2 X^{8b+4}+(\text{terms of degree $>8b+4$ in $X$})]\\
&=(-1/2)[\mu X^{4a+4b+4}+(\text{terms of degree $>4a+4b+4$ in $X$})]\\
&\quad\pm[X^{4a+2}
-(3/8)\mu^2 X^{4a+8b+6}+(\text{terms of degree $>4a+8b+6$ in $X$})]\\
\endalign
$$
where for the third equality we are using the fact that 
$\text{ord}_Xv(X)=a+b+1$ and hence
$$
v(X^4)=\mu X^{4a+4b+4}+(\text{terms of degree $>4a+4b+4$ in $X$})
\text{ with }0\ne\mu\in k
$$
and for the fourth equality we are using the Binomial Theorem for
exponent $1/2$. It follows that
$$
f^{(2)}(X^4,Y)=[Y-y_1(X)][Y-y_2(X)]
\text{ where }y_1\in k((X))\text{ and }y_2(X)\in k((X))
$$
are such that
$$
y_1(X)=X^{4a+2}-(\mu/2)X^{4a+4b+4}
+(\text{terms of degree $>4a+4b+4$ in $X$})
$$
and
$$
y_2(X)=-X^{4a+2}-(\mu/2)X^{4a+4b+4}
+(\text{terms of degree $>4a+4b+4$ in $X$}).
$$
We also have
$$
f^{(1)}(X^4,Y)=Y-v(X)\text{ where }v(X)\in k((X))\text{ with ord}_X=4a+4b+4.
$$
Finally by (TR4) of Section 6 we have
$$
f(X^4,Y)=\prod_{1\le j\le 4}[Y-z_j(X)]
$$
with
$$
z_j(X)=(\iota^jX)^{4a+2}+\frac{1}{2}(\iota^jX)^{4a+2b+3}
+(\text{terms of degree $>4a+2b+3$ in $X$})
$$
where $\iota$ is a primitive $4$-th root of $1$ in $k$. 
By the above expressions for the roots of $f$ and $f^{(1)}$ we get
$$
\text{int}(f,f^{(1)})=4a+2\text{ and noc}(f,f^{(1)})=(2a+1)/2.
$$
Likewise, by the above expressions for the roots of $f$ and $f^{(2)}$ we get
$$
\text{int}(f,f^{(2)})=8a+2b+5\text{ and noc}(f,f^{(2)})=(4a+2b+3)/4.
$$
It follows that
$$
\Omega'_{B_i}(f_Y)=f^{(i)}\text{ with deg}_Yf^{(i)}=D'(B_i)
\text{ for }1\le i\le 2
$$
and this completes the verification.

\centerline{}

\noindent
{\bf Example (DF5).} Finally, let us turn to the case of 
$F\in\widehat R^{\natural}$ having two factors, i.e., such that
$F=ff'$ with
$f\in\ R^{\natural}$ and $f'\in R^{\natural}$. At the same time let us 
arrange matters so that $F$ is pure meromorphic, i.e., 
$$
F(X,Y)=\Phi(X^{-1},Y)\text{ with }\Phi(X,Y)\in k[X,Y].
$$
To do this, in Example (TR5) of Section 6, let us take
$n>1$ with $b=0$ and $a=-1$, and
$$
\Phi(X,Y)=Y^{n+2}+\kappa' XY^n+\widehat\kappa Y+\kappa
+\sum_{3\le i\le n}\kappa_i Y^{n+2-i}\in k[X,Y]
$$
with
$$
0\ne\kappa'\in k\text{ and }0\ne\widehat\kappa\in k\text{ and }
0\ne\kappa\in k\text{ and }\kappa_i\in k\text{ for }3\le i\le n.
$$
As explained in (TR5) of Section 6, we then have
$$
F(X,Y)=f(X,Y)f'(X,Y) \text{ with }f(X,Y)\ne f'(X,Y)
$$
where
$$
f(X,Y)=Y^n+\sum_{1\le i\le n}w_i(X)Y^{n-i}\in R^{\natural}
\text{ and }w_i(X)\in k((X))
$$
with
$$
\text{ord}_Xw_i(X)>ie/n\text{ for $1\le i\le n-1$ and ord}_Xw_n(X)=e=1
$$
and
$$
f'(X,Y)=Y^2+\sum_{1\le i\le 2}w'_i(X)Y^{2-i}\in R^{\natural}
\text{ and }w'_i(X)\in k((X))
$$
with
$$
\text{ord}_Xw'_1(X)>e'/2\text{ and ord}_Xw'_2(X)=e'=-1
$$
and $0\ne\kappa'\in k$ and $0\ne\kappa/\kappa'\in k$ are the coefficients
of $X^{e'}$ and $X^e$ in $w'_2(X)$ and $w_n(X)$ respectively.
As explained in (TR5) of Section 6, we also have
$h(T(F))=2$ with $l_0(T(F))=-\infty$ and 
$$
l_1(T(F))=-1/2\text{ and }l_2(T(F))=1/n
$$
and upon letting
$$
\cases
B_0\in\widehat R^{\flat}\text{ with }
\sigma(B_0)=\{f,f'\}\text{ and }\lambda(B_0)=l_0(T(F)),\\
\text{and }B_1\in\widehat R^{\flat}\text{ with }
\sigma(B_1)=\{f,f'\}\text{ and }\lambda(B_1)=l_1(T(F)),\\
\text{and }B_2\in\widehat R^{\flat}\text{ with }
\sigma(B_2)=\{f\}\text{ and }\lambda(B_2)=l_2(T(F)),\\
\text{and }B'_2\in\widehat R^{\flat}\text{ with }
\sigma(B'_2)=\{f'\}\text{ and }\lambda(B'_2)=l_2(T(F)),\\
\endcases
$$ 
we have
$$
T(F)=\{B_0,B_1,B_2,B'_2\}
$$
with $D'(B_0)=0$ and 
$$
D'(B_1)=2\text{ and }D'(B_2)=n-1\text{ and }D'(B'_2)=0.
$$

Now
$$
F_Y=(n+2)Y^{n+1}+n\kappa' X^{-1}Y^{n-1}+\widehat\kappa
+\sum_{3\le i\le n}(n+2-i)\kappa_i Y^{n+1-i}
$$
and hence by (TR5) of Section 6 we see that $F_Y=(n+2)f^{(1)}f^{(2)}$
with $f^{(1)}\ne f^{(2)}$ where
$$
f^{(2)}(X,Y)
=Y^{n-1}+\sum_{1\le i\le n-1}v_i(X)Y^{n-1-i}\in R^{\natural}
\text{ and }v_i(X)\in k((X))
$$
with
$$
\text{ord}_Xv_i(X)>ie/n\text{ for }1\le i\le n-2 
\text{ and ord}_Xv_{n-1}(X)=e=1
$$
and
$$
f^{(1)}(X,Y)=Y^2+\sum_{1\le i\le 2}v'_i(X)Y^{2-i}\in R^{\natural}
\text{ and }v'_i(X)\in k((X))
$$
with
$$
\text{ord}_Xv'_1(X)>e'/2\text{ and ord}_Xv'_2(X)=e'=-1
$$
and $0\ne n\kappa'/(n+2)\in k$ and $0\ne\widehat\kappa/(n\kappa')\in k$ are 
the coefficients of $X^{e'}$ and $X^e$ in $v'_2(X)$ and 
$v_{n-1}(X)$ respectively.
In view of (TR3) of Section 6 we see that
$$
f(X^n,Y)=\prod_{1\le j\le n}[Y-z_j(X)]
$$
where $z_j(X)\in k((X))$ is such that
$$
z_j(X)=\omega^j \kappa^* X+(\text{terms of degree $>1$ in $X$})
$$
where $\omega$ is a primitive $n$-th root of $1$ in $k$, and $\kappa^*$
is an $n$-th root of $-\kappa/\kappa'$ in $k$, and
$$
f'(X^2,Y)=\prod_{1\le j\le 2}[Y-z'_j(X)]
$$
where $z'_j(X)\in k((X))$ is such that
$$
z'_j(X)=(-1)^j \kappa'{}^* X^{-1}+(\text{terms of degree $>-1$ in $X$})
$$
where $\kappa'{}^*$ is a square root of $-\kappa'$ in $k$.
In view of (TR3) of Section 6 we also see that
$$
f^{(2)}(X^{n-1},Y)=\prod_{1\le j\le n-1}[Y-y_j(X)]
$$
where $y_j(X)\in k((X))$ is such that
$$
y_j(X)=\widehat\omega^j\widehat\kappa^* X
+(\text{terms of degree $>1$ in $X$})
$$
where $\widehat\omega$ is a primitive $(n-1)$-th root of $1$ in $k$, 
and $\widehat\kappa^*$
is an $(n-1)$-th root of $-\widehat\kappa/(n\kappa')$ in $k$, and
$$
f^{(1)}(X^2,Y)=\prod_{1\le j\le 2}[Y-y'_j(X)]
$$
where $y'_j(X)\in k((X))$ is such that
$$
y'_j(X)=(-1)^j\widehat\kappa'{}^* X^{-1}
+(\text{terms of degree $>-1$ in $X$})
$$
where $\widehat\kappa'{}^*$ is a square root of $-n\kappa'/(n+2)$ in $k$.
By the above expressions of the roots of $f,f',f^{(1)},f^{(2)}$ we get
$$
\cases
\text{int}(f,f^{(1)})=-n\text{ and noc}(f,f^{(1)})=-1/2,\text{ and }\\
\text{int}(f',f^{(1)})=-2\text{ and noc}(f',f^{(1)})=-1/2,\text{ and }\\
\text{int}(f,f^{(2)})=(n-1)\text{ and noc}(f,f^{(2)})=1/n,\text{ and }\\
\text{int}(f',f^{(2)})=-(n-1)\text{ and noc}(f',f^{(2)})=-1/2,
\endcases
$$
and hence
$$
\Omega'_{B_i}(F_Y)=f^{(i)}\text{ with deg}_Yf^{(i)}=D'(B_i)
\text{ for }1\le i\le 2
$$
which verifies Theorem (DF1) in the present situation.

To get an example of $\widehat F\in\widehat R^{\natural}$ having three
factors, we take 
$\widehat F(X,Y)=\widehat\Phi(X^{-1},Y)$ with
$\widehat\Phi(X,Y)=\Phi(X,Y)-\Phi(0,0)\in k[X,Y]$. Then by
(TR5) of Section 6 we get 
$\widehat F=\widehat f\widehat f'\widehat f''$
with $\widehat f''\ne\widehat f\ne\widehat f'\ne\widehat f''$
where $\widehat f''(X,Y)=Y\in R^{\natural}$ and
$$
\widehat f(X,Y)
=Y^{n-1}+\sum_{1\le i\le n-1}\widehat w_i(X)Y^{n-1-i}\in R^{\natural}
\text{ and }\widehat w_i(X)\in k((X))
$$
with
$$
\text{ord}_X\widehat w_i(X)>ie/n\text{ for }1\le i\le n-2 
\text{ and ord}_X\widehat w_{n-1}(X)=e=1
$$
and
$$
\widehat f'(X,Y)
=Y^2+\sum_{1\le i\le 2}\widehat w'_i(X)Y^{2-i}\in R^{\natural}
\text{ and }\widehat w'_i(X)\in k((X))
$$
with
$$
\text{ord}_X\widehat w'_1(X)>e'/2\text{ and ord}_X\widehat w'_2(X)=e'=-1
$$
and $0\ne \kappa'\in k$ and $0\ne\widehat\kappa/\kappa'\in k$ are 
the coefficients of $X^{e'}$ and $X^e$ in $\widehat w'_2(X)$ and 
$\widehat w_{n-1}(X)$ respectively.
In view of (TR3) of Section 6 we see that
$$
\widehat f(X^{n-1},Y)=\prod_{1\le j\le n-1}[Y-\widehat z_j(X)]
$$
where $\widehat z_j(X)\in k((X))$ is such that
$$
\widehat z_j(X)=\widehat\omega^j\widehat\kappa^* X
+(\text{terms of degree $>1$ in $X$})
$$
where $\widehat\omega$ is a primitive $(n-1)$-th root of $1$ in $k$, 
and $\widehat\kappa^*$
is an $(n-1)$-th root of $-\widehat\kappa/\kappa'$ in $k$, and
$$
\widehat f'(X^2,Y)=\prod_{1\le j\le 2}[Y-\widehat z'_j(X)]
$$
where $\widehat z'_j(X)\in k((X))$ is such that
$$
\widehat z'_j(X)=(-1)^j\widehat\kappa'{}^* X^{-1}
+(\text{terms of degree $>-1$ in $X$})
$$
where $\widehat\kappa'{}^*$ is a square root of $-\kappa'$ in $k$. 
By the above expressions of the roots of $\widehat f$ and $\widehat f'$
it follows that
$h(T(\widehat F))=2$ with $l_0(T(\widehat F))=-\infty$ and 
$$
l_1(T(\widehat F))=-1/2\text{ and }l_2(T(\widehat F))=1/(n-1)
$$
and upon letting
$$
\cases
\widehat B_0\in\widehat R^{\flat}
\text{ with }\sigma(\widehat B_0)=\{\widehat f,\widehat f',\widehat f''\}
\text{ and }\lambda(\widehat B_0)=l_0(T(\widehat F)),\\
\text{and }\widehat B_1\in\widehat R^{\flat}
\text{ with }\sigma(\widehat B_1)=\{\widehat f,\widehat f',\widehat f''\}
\text{ and }\lambda(\widehat B_1)=l_1(T(\widehat F)),\\
\text{and }\widehat B_2\in\widehat R^{\flat}
\text{ with }\sigma(\widehat B_2)=\{\widehat f,\widehat f''\}
\text{ and }\lambda(\widehat B_2)=l_2(T(\widehat F)),\\
\text{and }\widehat B'_2\in\widehat R^{\flat}
\text{ with }\sigma(\widehat B'_2)=\{\widehat f'\}
\text{ and }\lambda(\widehat B'_2)=l_2(T(\widehat F)),\\
\endcases
$$ 
we have
$$
T(\widehat F)=\{\widehat B_0,\widehat B_1,\widehat B_2,\widehat B'_2\}
$$
with $D'(\widehat B_0)=0$ and 
$$
D'(\widehat B_1)=2
\text{ and }D'(\widehat B_2)=n-1\text{ and }D'(\widehat B'_2)=0.
$$
Now $\widehat F_Y=F_Y=(n+2)f^{(1)}f^{(2)}$, and by the above expressions 
of the roots of $\widehat f,\widehat f',f^{(1)},f^{(2)}$ we get
$$
\cases
\text{int}(\widehat f,f^{(1)})=-(n-1)
\text{ and noc}(\widehat f,f^{(1)})=-1/2,\text{ and }\\
\text{int}(\widehat f',f^{(1)})=-2
\text{ and noc}(\widehat f',f^{(1)})=-1/2,\text{ and }\\
\text{int}(\widehat f'',f^{(1)})=-1
\text{ and noc}(\widehat f'',f^{(1)})=-1/2,\text{ and }\\
\text{int}(\widehat f,f^{(2)})=(n-1)
\text{ and noc}(\widehat f,f^{(2)})=1/(n-1),\text{ and }\\
\text{int}(\widehat f',f^{(2)})=-(n-1)
\text{ and noc}(\widehat f',f^{(2)})=-1/2,\text{ and }\\
\text{int}(\widehat f'',f^{(2)})=1
\text{ and noc}(\widehat f'',f^{(2)})=1/(n-1),
\endcases
$$
and hence
$$
\Omega'_{\widehat B_i}(\widehat F_Y)
=f^{(i)}\text{ with deg}_Yf^{(i)}=D'(\widehat B_i)
\text{ for }1\le i\le 2
$$
which again verifies Theorem (DF1) in the present situation.

\centerline{}

\centerline{{\bf Section 9: Factorization of the Jacobian}}

If $T=T(FG)$ where $F\in R\setminus k((X))$ and $0\ne G\in R$, 
then for every $B\in T$ with $\Omega_B(G)=1$
we clearly have $R^*(F,B)=R^*(B)$. Therefore by (SP76), (SP77),
(SP79) and (SP80) we get the following Jacobian Factorization Theorem.

\proclaim{Theorem (JF1)}
Let $T=T(FG)$ where $F\in R\setminus k((X))$ is devoid of multiple 
factors in $R$, and $0\ne G\in R$. Then we have the following. 

(JF1.1) If $B\in T$ is such that $\Omega_B(G)=1$ and $S(G,B)\ne 0$
then we have
$$
\text{deg}_Y\Omega_B(J(F,G))
=\text{deg}_Y\Omega_B((FG)_Y)
=\text{deg}_Y\Omega_B(F_Y)=D''(B)
$$
and
$$
\text{deg}_Y\Omega'_B(J(F,G))
=\text{deg}_Y\Omega'_B((FG)_Y)
=\text{deg}_Y\Omega'_B(F_Y)=D'(B)
$$
and
$$
\Omega_B(J(F,G))
=\Omega'_B(J(F,G))\prod_{B'\in\pi(T,B)}\Omega'_{B'}(J(F,G))
$$
where for every $B'\in\pi(T,B)$ we have
$\Omega_{B'}(G)=1$ and $S(G,B')\ne 0$ and
$$
\text{deg}_Y\Omega_{B'}(J(F,G))
=\text{deg}_Y\Omega_{B'}((FG)_Y)
=\text{deg}_Y\Omega_{B'}(F_Y)=D''(B')
$$
and
$$
\text{deg}_Y\Omega'_{B'}(J(F,G))
=\text{deg}_Y\Omega'_{B'}((FG)_Y)
=\text{deg}_Y\Omega'_{B'}(F_Y)=D'(B').
$$

(JF1.2) If $B\in T$ is such that $\Omega_B(G)=1$
then we have
$$
\text{deg}_Y\Omega_B((FG)_Y)
=\text{deg}_Y\Omega_B(F_Y)=D''(B)
$$
and
$$
\text{deg}_Y\Omega'_B((FG)_Y)
=\text{deg}_Y\Omega'_B(F_Y)=D'(B)
$$
and
$$
\Omega_B((FG)_Y)
=\Omega'_B((FG)_Y)\prod_{B'\in\pi(T,B)}\Omega'_{B'}((FG)_Y)
$$
where for every $B'\in\pi(T,B)$ we have
$\Omega_{B'}(G)=1$ and 
$$
\text{deg}_Y\Omega_{B'}((FG)_Y)
=\text{deg}_Y\Omega_{B'}(F_Y)=D''(B')
$$
and
$$
\text{deg}_Y\Omega'_{B'}((FG)_Y)
=\text{deg}_Y\Omega'_{B'}(F_Y)=D'(B').
$$
\endproclaim

\centerline{}

\noindent
{\bf Remark (JF2).}
The Jacobian Factorization (JF1.1) was based on (SP80), and it invoked the
Jacobian Estimates (JE1) to (JE3) of Section 5 only in the special case
when $\text{deg}_Y\text{minco}_X[z,V,W](G)=0$. Elsewhere we shall discuss
a more refined Jacobian Factorization based on (SP30) by invoking the
general case of (JE1) to (JE3).

\centerline{}

\noindent
{\bf Example (JF3).} 
Now let us illustrate Theorem (JF1.1) by the example
$$
F=F(X,Y)=Y^n+X^e\in R^{\natural}
\text{ where }0\ne e\in\Bbb Z\text{ with GCD}(n,e)=1
$$
considered in (DF3) of Section 8.
For $G=X$ we have $J(F,G)=F_Y$ and we are reduced to (DF3).

\centerline{}

\noindent
{\bf Example (JF4).} 
Next let us illustrate Theorem (JF1.1) by the example
$$
F=F(X,Y)=(Y^2-X^{2a+1})^2-X^{3a+b+2}Y\in R^{\natural}
\text{ with }a\in\Bbb Z\text{ and }0\le b\in\Bbb Z
$$
considered in (DF4) of Section 8.
Again, for $G=X$ we have $J(F,G)=F_Y$ and we are reduced to (DF4). 
Moreover, for 
$$
\widehat G=\widehat G(X,Y)=Y\in R^{\natural}
$$
we have 
$$
J(F,\widehat G)=F_X=-(4a+2)X^{2a}\widehat F
$$ 
with
$$
\widehat F=\widehat F(X,Y)=Y^2+(3a+b+2)(4a+2)^{-1}X^{a+b+1}-X^{2a+1}.
$$
By (TR3) and (TR4) of Section 6, it follows that
$h(T(F\widehat G))=2$ with $l_0(T(F\widehat G))=-\infty$ and 
$$
l_1(T(F\widehat G))=(2a+1)/2\text{ and }l_2(T(F\widehat G))=(4a+2b+3)/4
$$
and upon letting
$$
\cases
\widehat B_0\in\widehat R^{\flat}
\text{ with }\sigma(\widehat B_0)=\{F,\widehat G\}
\text{ and }\lambda(\widehat B_0)=l_0(T(F\widehat G)),\\
\text{and }\widehat B_1\in\widehat R^{\flat}
\text{ with }\sigma(\widehat B_1)=\{F,\widehat G\}
\text{ and }\lambda(\widehat B_1)=l_1(T(F\widehat G)),\\
\text{and }\widehat B_2\in\widehat R^{\flat}
\text{ with }\sigma(\widehat B_2)=\{F\}
\text{ and }\lambda(\widehat B_2)=l_2(T(F\widehat G)),\\
\text{and }\widehat B'_2\in\widehat R^{\flat}
\text{ with }\sigma(\widehat B'_2)=\{\widehat G\}
\text{ and }\lambda(\widehat B'_2)=l_2(T(F\widehat G)),\\
\endcases
$$ 
we have
$$
T(F\widehat G)=\{\widehat B_0,\widehat B_1,\widehat B_2,\widehat B'_2\}
$$
with $D'(\widehat B_0)=0$ and 
$$
D'(\widehat B_1)=2
\text{ and }D'(\widehat B_2)=2\text{ and }D'(\widehat B'_2)=0.
$$
By (TR3) and (TR4) we also see that $\widehat F\in R^{\natural}$ with
$$
\text{noc}(\widehat G,\widehat F)=(2a+1)/4\text{ and }
\text{noc}(F,\widehat F)=(4a+2b+3)/4
$$
and hence
$$
\Omega'_B(J(F,\widehat G))=\cases
\widehat F&\text{if }B=\widehat B_2\\
1&\text{if $B=\widehat B_0$ or $B=\widehat B_1$ or }B=\widehat B'_2
\endcases
$$
in accordance with Theorem (JF1.1).
Likewise, for 
$$
\widetilde G=\widetilde G(X,Y)=Y^2-X^{2a+1}\in R^{\natural}
$$ 
we have
$$
J(F,\widetilde G)=J(-X^{3a+b+2}Y,\widetilde G)=\cases
-(6a+2b+4)X^{3a+b+1}\widetilde F&\text{if }3a+b+2\ne 0\\
-(2a+1)X^{2a}\widetilde F&\text{ if }3a+b+2=0\\
\endcases
$$
with
$$
\widetilde F=\widetilde F(X,Y)=\cases
Y^2+(2a+1)(6a+2b+4)^{-1}X^{2a+1}&\text{if }3a+b+2\ne 0\\ 
1&\text{if }3a+b+2=0\\ 
\endcases
$$
By (TR3) and (TR4) of Section 6, it follows that
$h(T(F\widetilde G))=2$ with $l_0(T(F\widetilde G))=-\infty$ and 
$$
l_1(T(F\widetilde G))=(2a+1)/2\text{ and }l_2(T(F\widetilde G))=(4a+2b+3)/4
$$
and upon letting
$$
\cases
\widetilde B_0\in\widehat R^{\flat}
\text{ with }\sigma(\widetilde B_0)=\{F,\widetilde G\}
\text{ and }\lambda(\widetilde B_0)=l_0(T(F\widetilde G)),\\
\text{and }\widetilde B_1\in\widehat R^{\flat}
\text{ with }\sigma(\widetilde B_1)=\{F,\widetilde G\}
\text{ and }\lambda(\widetilde B_1)=l_1(T(F\widetilde G)),\\
\text{and }\widetilde B_2\in\widehat R^{\flat}
\text{ with }\sigma(\widetilde B_2)=\{F,\widetilde G\}
\text{ and }\lambda(\widetilde B_2)=l_2(T(F\widetilde G)),\\
\endcases
$$ 
we have
$$
T(F\widetilde G)=\{\widetilde B_0,\widetilde B_1,\widetilde B_2\}
$$
with $D'(\widetilde B_0)=0$ and 
$$
D'(\widetilde B_1)=3\text{ and }D'(\widetilde B_2)=4.
$$
Thus the stem of every bud of $T(F\widetilde G)$ contains $F$ as well as
$\widetilde G$, and hence Theorem (JF1.1) does not predict any factors of
$J(F,\widetilde G)$. This is quite satisfactory when $3a+b+2=0$ because
then $J(F,\widetilde G)=-(2a+1)X^{2a+1}$ and so 
$J(F,\widetilde G)$ has no factor involving $Y$. A particularly
interesting case of $3a+2b+2=0$ is the pure meromorphic
case when $(a,b)=(-1,1)$.
In that case, as noted in (TR4) of Section 6,
we have $F(X,Y)=\Phi(X^{-1},Y)$ where $\Phi(X,Y)\in k[X,Y]$ is the
{\bf variable} $\Phi(X,Y)=(Y^2-X)^2-Y$; indeed, then
$k[X,Y]=k[\Phi,\Psi]$ where $\Psi(X,Y)=Y^2-X$ with
$\widetilde G(X,Y)=\Psi(X^{-1},Y)$.

\centerline{}

\noindent
{\bf Example (JF5).} 
Finally let us illustrate Theorem (JF1.1) by the example
$$
F=F(X,Y)=\Phi(X^{-1},Y)
$$
where
$$
\Phi(X,Y)=Y^{n+2}+\kappa' XY^n+\widehat\kappa Y+\kappa
+\sum_{3\le i\le n}\kappa_i Y^{n+2-i}\in k[X,Y]
$$
with $n>1$ and
$$
0\ne\kappa'\in k\text{ and }0\ne\widehat\kappa\in k\text{ and }
\kappa\in k\text{ and }\kappa_i\in k\text{ for }3\le i\le n
$$
considered in (DF5) of Section 8 (where we used the notation $\widehat\Phi$
and $\widehat F$ for the special case of $\kappa=0$).
Once again, for $G=-X$ we have $J(F,G)=F_Y$ and we are reduced to (DF5).
Note that the affine plane curve $\Phi=0$ is nonsingular (at finite distance)
for every $\kappa$; equivalently, $\Phi_X$ and $\Phi_Y$ have no solution in
the affine plane $k\times k$.
Moreover, $\Phi$ is irreducible for every $\kappa\ne 0$, but
reducible for $\kappa=0$. 
In (DF5) we have shown that $\Phi$ has two or three places at $\infty$,
i.e., $F$ has two or three factors in $R^{\natural}$, according as
$\kappa\ne 0$ or $\kappa=0$.
Further interest in this nice family of bivariate polynomials $\Phi$
lies in the fact that it provides a convenient
testing ground for the trivariate jacobian conjecture. 
To elucidate this, given any $H_1\in k[X_1,\dots,X_r]$ where
$r$ is any positive integer, let us say that $H_1$ is a {\bf variable in}
$k[X_1,\dots,X_r]$ to mean that $k[X_1,\dots,X_r]=k[H_1,\dots,H_r]$ for
some $H_2,\dots,H_r$ in $k[X_1,\dots,X_r]$, and 
let us say that $H_1$ is a {\bf weak variable in}
$k[X_1,\dots,X_r]$ to mean that $0\ne J(H_1,\dots,H_r)\in k$
for some $H_2,\dots,H_r$ in $k[X_1,\dots,X_r]$ where $J(H_1\dots,H_r)$ is
the jacobian of $H_1,\dots,H_r$ with respect to $X_1,\dots,X_r$.
The reducibility of $\Phi$ when $\kappa=0$ shows that $\Phi$ is not a
variable in $k[X,Y]$. 
The reducibility of $\Phi$ when $\kappa=0$ also shows that $\Phi$ is not a
variable in $k[X,Y,Z]$. 
It can be shown that $\Phi$ is not a weak variable in $k[X,Y]$, 
at least when $n+1$ is a prime number.
However, as was pointed out to us by Ignacio Luengo,
it is not known whether $\Phi$ is or is not a weak variable in 
$k[X,Y,Z]$, even when $n=2$.
To see that $\Phi$ is not a weak variable in $k[X,Y]$, first note that
by the automorphism
$(X,Y)\mapsto((X-Y^2)/\kappa',Y)$ we can send $\Phi$ to the polynomial
$$
XY^n+\widehat\kappa Y+\kappa
+\sum_{3\le i\le n}\kappa_i Y^{n+2-i}\in k[X,Y]
$$
whose degree is $n+1$ and whose degree form $XY^n$ has two coprime
factors. On the other hand, 
it can easily be shown that if $H\in k[X,Y]$ is a weak variable
in $k[X,Y]$ of prime degree then its degree form must be a power of a
linear form. 

\centerline{}

\centerline{{\bf REFERENCES}}

\widestnumber\key{AOS}

\ref \key Ab
\by S. S. Abhyankar
\paper On the semigroup of a meromorphic curve, Part I
\jour Proceedings of the International Symposium of Algebraic Geometry, Kyoto,
\yr 1977
\pages 240-414
\endref

\ref \key De
\by F. Delgado de la Mata
\paper A factorization theorem for the polar of a curve with two branches
\jour Compositio Mathematica
\vol 92
\yr 1994
\pages  327-375
\endref

\ref \key KL
\by T.C. Kuo and Y.C. Lu
\paper On analytic function germs of two complex variables
\jour Topology
\vol 16
\yr 1977
\pages 299-310
\endref

\ref \key Me
\by M. Merle
\paper Invariants polaires des courbes planes
\jour Inventiones Mathematicae
\vol 41
\yr 1977
\pages  299-310
\endref

\enddocument